%% file: main_MSOM.tex
\documentclass[msom]{oo}

\OneAndAHalfSpacedXI 

\usepackage{natbib}
 \bibpunct[, ]{(}{)}{,}{a}{}{,}%
 \def\bibfont{\small}%
\TheoremsNumberedThrough     
\EquationsNumberedThrough 
\input{Definitions} 

\begin{document}
\RUNAUTHOR{Sereshti et al.}
\RUNTITLE{Stochastic Dynamic Lot-sizing with Supplier-Driven Substitution and Service Level Constraints}
\TITLE{Stochastic Dynamic Lot-sizing with Supplier-Driven Substitution and Service Level Constraints}
\ARTICLEAUTHORS{%
\AUTHOR{Narges Sereshti}
\AFF{Department of Decision Sciences, HEC Montr\'{e}al, Montr\'{e}al, Qu\'{e}bec H3T 2A7, Canada, \EMAIL{narges.sereshti@hec.ca}}
\AUTHOR{Merve Bodur}
\AFF{Department of Mechanical and Industrial Engineering, University of Toronto, Toronto, Ontario M5S 3G8, Canada,
\EMAIL{bodur@mie.utoronto.ca}}
\AUTHOR{James R. Luedtke}
\AFF{Department of Industrial and Systems Engineering, University of Wisconsin, Madison, Wisconsin 53706, \EMAIL{jim.luedtke@wisc.edu}}
}
\ABSTRACT{%
We consider a multi-stage stochastic lot-sizing problem with service level constraints and supplier-driven product substitution. A firm has multiple products and it has the option to meet demand from substitutable products at a cost. Considering the uncertainty in future demands, the firm wishes to make ordering decisions in every period such that the probability that all demands can be met in the next period meets or exceeds a minimum service level. 
We propose a rolling-horizon policy in which a two-stage joint chance-constrained stochastic program is solved to make decisions in each time period. We demonstrate how to effectively solve this formulation. In addition, we propose two policies based on deterministic approximations. 
We demonstrate that the proposed chance-constraint policy can achieve the service levels more reliably and at a lower cost. We also explore the value of product substitution in this model, demonstrating that the substitution option allows achieving service levels while reducing costs by 7\% to 25\% in our experiments, and that the majority of the benefit can be obtained with limited levels of substitution allowed.
}%
\HISTORY{}
\maketitle
\section{Introduction}
The basic lot-sizing problem is a multi-period production planning problem that considers the trade-off between setup costs and inventory holding costs and defines the optimal timing and quantity of production to minimize the total cost. 
When demand is uncertain, which is inevitable in real-world applications, the decision-maker needs to determine the production policy to minimize the expected cost over the distribution of demand outcomes. Demand uncertainty leads to the possibility of stock-outs (i.e., demand exceeds available inventory) and a key challenge then is to limit the frequency of such undesirable events. A common approach to deal with this challenge is to assume customer demands can be backlogged (i.e., met in a period later than when it arrived) and assign a cost per period that it is backlogged. Then, the incurred backlog cost needs to capture the costs associated with both tangible and intangible effects which may be difficult to estimate. 
In contrast, in this work we study the stochastic lot-sizing problem with an $\alpha$ service level constraint which instead requires that there is no stock-out in each period with probability at least $\alpha$. 

The standard strategy for limiting stock-outs is to hold more product in inventory, which leads to a trade-off between inventory holding costs and service level. In some cases, when a firm is managing inventory and ordering decisions of multiple products the firm has the option to substitute one product for another to avoid a stock-out. This type of substitution is known as supplier-driven substitution and provides another mechanism to avoid stock-outs which can be interpreted as a risk-pooling strategy for handling uncertain demand~\citep{shin2015classification}.  Supplier-driven substitution has practical relevance in the electronics and steel industries where it is possible to substitute a lower-grade product with a higher-grade one~\citep{lang2010efficient}. 

To explore the potential benefits of supplier-driven substitution, we study the stochastic lot-sizing problem with substitution and joint service level constraint over multiple products. A joint service level constraint 
ensures that no products have a stock-out to exceed the target $\alpha$ in each period. A joint service level is necessary in our model because, given on-hand inventory and observed costumer demands, we jointly determine what substitutions should be made in order to avoid a stock-out. Such joint determination links the stock-out event of different products together, making it impossible to separately control the probability of each individual product having a stock-out. Note that, if the substitution policy is fixed (e.g., one always substitutes product 1 for product 2 if there is a shortfall in product 2, etc.) then it would be possible to constrain individual product service levels. We do not pursue this option, as we prefer to allow substitution decisions to be flexibly optimized in each period given the available inventory and current product demands. 


We consider an infinite-horizon problem in which the firm sequentially makes setup, production, and substitution decisions based on the current state of the system, reflected as the amount of available inventory and backlog amounts of each product. We follow the ``dynamic" strategy \citep{bookbinder1988strategies} in which decisions are made throughout the planning horizon in response to the latest observed information. As the infinite-horizon problem is computationally intractable, we propose to solve a finite-horizon problem and apply it in a rolling-horizon fashion. The aim is to propose decision policies that make the current period decisions by solving a finite-horizon problem that looks ahead a certain number of periods. 
Ideally, this finite-horizon problem would take the form of a multi-stage stochastic program that considers all possible sample paths over the horizon and anticipates the future optimal decisions. While we begin by formulating this ideal model, it is also computationally intractable.
Thus, we propose to solve approximations of this problem to drive our decision policy. The first approximation is purely deterministic, as commonly employed in practice for various application domains, and hence is not able to explicitly consider the service level constraint. We then propose a chance-constrained approximation that considers scenarios of possible joint demands in the next stage and hence is able to enforce that the chosen decisions satisfy the service level constraint. This two-stage model contains a joint chance constraint, for which we apply results in \citep{luedtke2014branch} to derive an efficient branch-and-cut (B\&C) algorithm. This model differs from standard two-stage approximations of multi-stage stochastic programming models in that it considers a distribution of scenarios of product demands in the immediate next stage, but for stages beyond that it merges these approximations back into a deterministic approximation, which is done to improve tractability of the model. 

We use simulation to evaluate our proposed policies in a steady-state system and demonstrate that over a range of problem characteristics the policy driven by our proposed chance-constrained model respects the service level targets more reliably and at a lower cost than the policies driven by solving deterministic models. We also explore the value of product substitution and find that using substitution achieves the target service levels at significantly reduced costs compared to without substitution and that the majority of the benefits can be obtained even when limiting substitution to be between products of the most similar quality.

We summarize our main contributions as follows. 
\begin{itemize}
\item We study an infinite-horizon multi-stage lot-sizing problem with substitution and joint service level constraints, which to the best of our knowledge is new to the literature.

\item We propose rolling-horizon policies based on solving finite-horizon deterministic and two-stage chance-constrained optimization models to make decisions in each period.

\item We describe a branch-and-cut algorithm to solve the two-stage chance-constrained optimization model and demonstrate its computational efficiency.

\item We conduct a simulation study that demonstrates the value of the chance-constrained optimization driven policy and the value of supplier-driven substitution. We also provide insights obtained through sensitivity analysis on important parameters of the problem, such as when substitution is most valuable. 

\end{itemize}

The rest of the paper is organized as follows. In Section \ref{sec:litrev}, we review the related literature. In Section \ref{sec:form}, we define the problem and the dynamics of decisions in the system and provide a dynamic programming formulation for the finite-horizon multi-stage stochastic program. In Section \ref{sec:approx}, we describe the optimization models (approximations of the model from Section 3) that we propose to use to make decisions and present the B\&C algorithm to solve the chance-constrained model.
In Section \ref{sec:comp}, we present results from the computational experiments, including policy comparison and insights about the value of substitution. We provide concluding remarks in Section \ref{sec:conc}.

\section{Literature review}
\label{sec:litrev}
The related literature to this work can be categorized in two streams. The first part is dedicated to the lot-sizing and inventory models with substitution in both deterministic and stochastic versions and the second part is dedicated to the stochastic lot-sizing problem with joint service level. In what follows, we review the related works, whose main characteristics are summarized in Table~\ref{tab:litRev}. To the best of our knowledge, no research has investigated the stochastic lot-sizing problem with substitution and joint service levels.


\input{LiteratureTable}

\subsection{Lot-sizing and inventory problems with substitution}
There are two types of substitution: customer-driven and supplier-driven~\citep{shin2015classification}. In customer-driven substitution, the customer decides which product to substitute~\citep{zeppetella2017optimal}, while in the supplier-driven (firm-driven) case, it is the supplier, firm, or the vendor who makes the substitution decisions~\citep{rao2004multi}. The substitution possibility is addressed in both deterministic and stochastic settings which are explained as follows.

\textbf{Deterministic models.}
\cite{hsu2005dynamic} study two versions of the dynamic uncapacitated lot-sizing problem with supplier-driven substitution, when there is a need for physical conversion before substitution, and when no conversion is needed. They propose a mixed-integer linear programming (MILP) model and solve it using a backward dynamic programming algorithm and an algorithm based on Silver-Meal heuristic. \cite{lang2010efficient} consider the uncapacitated lot-sizing problem with general substitution in which a specific class of demand can be satisfied by different products based on a substitution graph. They propose a MILP model along with some valid inequalities, also a plant location reformulation in which the amount of production for an item is broken down into different amounts based on the period where they are used to satisfy the demand. 

\textbf{Stochastic models.}
The majority of studies in stochastic inventory planning have considered customer-driven substitution. 
This is also known as ``stock-out substitution". \cite{akccaycategory} investigate a single-period inventory planning problem with substitutable products. Considering the stock-out substitution, they propose an optimization based method, which jointly defines the ordering decisions of each product, while satisfying a service level. 
\cite{nagarajan2008inventory} consider the inventory planning problem with customer-driven substitution. They propose an optimal policy for specific cases in terms of planning periods and the number of products, and a heuristic algorithm for the general form of the problem. In the inventory planning problem, there is no setup cost and they try to optimize the profit in the system which is equal to selling revenue minus the holding, substitution, and lost sales costs.

Our model considers supplier-driven substitution. \cite{bassok1999single} investigate the single-period inventory management problem with random demand and downward supplier-driven substitution in which a lower-grade item can be substituted with  ones with a higher-grade. This model is an extension of the newsvendor problem and there is no setup cost in case of ordering. 
The authors propose a profit maximization formulation and characterize the structure of the optimal policy for this single-period problem. They propose bounds on the optimal order amount and use them in an iterative algorithm to solve the model. \cite{rao2004multi} also consider a single-period problem with stochastic demand and downward substitution, and model it as a two-stage stochastic program. Their model considers the initial inventory and the ordering cost. They derive a deterministic equivalent formulation (extensive form) and propose two heuristic algorithms to solve this problem.

Another related research stream considers the possibility of having multiple graded output items from a single input item, which is known as ``co-production"~\citep{ng2012robust}. In these problems, there is a hierarchy in the grade of output items and it is possible to substitute a lower-grade item with the ones with higher-grade~\citep{bitran1992ordering}. \cite{hsu1999random} consider the single-period production system with random demand and random yields. 
They model the problem as a two-stage stochastic program which defines the production amount of a single item and the allocation of its output items to different demand classes. They propose  
two decomposition based methods, in which the subproblems are network flow problems. 
 \cite{bitran1992ordering} study an infinite-horizon, multi-item, multi-period co-production problem with deterministic demand and random yields. As solving this problem in an infinite-horizon is intractable, they propose two approximation approaches. 
The first approximation is based on a rolling-horizon implementation of the finite-horizon stochastic model, which is related to the overall approach we take. For the second approximation, they consider a simple heuristic based on the optimal allocation policy, in a multi-period setting. This heuristic includes two modules; a module to determine the production quantities, and a module to allocate produced items to the customers. This heuristic can be also applied in a rolling-horizon procedure. \cite{bitran1992deterministic} consider the same problem and propose deterministic near-optimal approximations within a fixed planning horizon. To adapt their model to the revealed information, they apply the proposed model using simple heuristics in a rolling planning horizon. \cite{bitran1994co} consider the co-production and random yield in a semiconductor industry and propose heuristic methods to solve their proposed model. 
 
\subsection{Stochastic lot-sizing problem and service level constraints}
Most of the research about stochastic lot-sizing problem with stochastic demands consider a scenario set or a scenario tree to represent the randomness in demand. Much of this research assumes backlogging has a cost that is included in the objective. \cite{haugen2001progressive} consider the multi-stage uncapacitated lot-sizing problem and propose a progressive hedging algorithm to solve their proposed model. \cite{guan2008polynomial} propose a dynamic programming algorithm for a similar model. Using the same algorithm, \cite{guan2011stochastic} studies the capacitated version of the problem with the possibility of backlogging. \cite{lulli2004branch} propose a branch-and-price algorithm for multi-stage stochastic integer programming and apply their general method to the stochastic batch-sizing problem. In this problem, they consider that the demand, production, inventory and setup costs are uncertain. The difference between this problem and the lot-sizing problem is that the production quantities are in batches and the production decisions are the integer-valued number of batches that will be produced. 
\cite{lulli2006heuristic} also proposed a heuristic scenario updating method for the stochastic batch-sizing problem.

Stochastic lot-sizing problems with service level constraints have been studied extensively~\citep{tempelmeier2007stochastic} and many types of service levels exist in the literature. One of the main service levels is the $\alpha$ service level which is an event-oriented service level, and imposes limits on the probability of a stock-out. This service level is represented as a chance-constraint and is usually defined for each period and product separately. \cite{bookbinder1988strategies} investigate stochastic lot-sizing problems with an $\alpha$ service level and propose three different strategies for this problem based on the timing of the setup and production decisions. These strategies are the \textit{static}, \textit{dynamic}, and \textit{static-dynamic} strategies. In the \textit{static} strategy, both the setup and production decisions are determined at the beginning of the planning horizon and they remain fixed when the demand is realized. In the \textit{dynamic} strategy, both the setup and production decisions are dynamically changed with the demand realizations throughout the planning horizon. The \textit{static-dynamic} strategy is between these two strategies in which the setups are fixed at the beginning of the planning horizon and the production decisions are updated when the demands are realized. In this work, we will follow the \textit{dynamic} strategy and all the decisions are updated dynamically with the demand realization. 

Some studies define the service level constraint jointly over periods in the planning horizon. \cite{liu2018polyhedral} consider the uncapacitated lot-sizing problem with a joint service level constraint. They study the polyhedral structure of the problem, and propose different valid inequalities and a reformulation of the  problem. \cite{zhangetal2022} extend this line of work with additional valid inequalities and formulations. \cite{jiang2017production} consider the same problem with and without pricing decisions. \cite{gicquel2018joint} investigate the capacitated version of the same problem. \cite{jiang2017production} and \cite{gicquel2018joint} use a sample average approximation method to solve their problems, which is a variation of the 
method proposed by \cite{luedtke2008sample} to solve models with chance-constraints using scenario sets. All of these studies consider single item models in which the joint service level is defined over all periods. Few studies consider the service level jointly over all the products. \cite{akccaycategory} adapt the Type II service level or ``fill rate" for each individual product and overall within a category of products in the customer-driven substitution model. This type of service level considers the expected value of backlog and it is not modeled as a chance constraint. \cite{sereshti2020value} study different types of aggregate service level for the lot-sizing problem which are defined over multiple products, but they do not consider substitution. In this work, we consider supplier-driven substitution and a joint service level that is defined over all products, but separately for each time period.

\section{Problem definition and formulation}
\label{sec:form}

We consider a stochastic lot-sizing problem with the possibility of supplier-driven substitution in an infinite time horizon which is discretized into planning periods. 
There are multiple types of products  $\KA= \{1,...,\Ka\}$ with random demand, and at each period, we need to make decisions about the production setups, production and substitution amounts, and accordingly define the inventory and backlog levels. There is a production lead time of one period, i.e., what is produced in the current period is available to meet demand in the next period or later. 
These decisions are made sequentially in each period based on the available inventory and backlog in the system in that period, such that a joint service level over all products is to be satisfied in the following period. This fits into the category of the  ``dynamic'' strategy that is defined for the stochastic lot-sizing problem~\citep{bookbinder1988strategies}. 

\begin{figure}[ht]
\begin{center}
\includegraphics[scale=0.6]{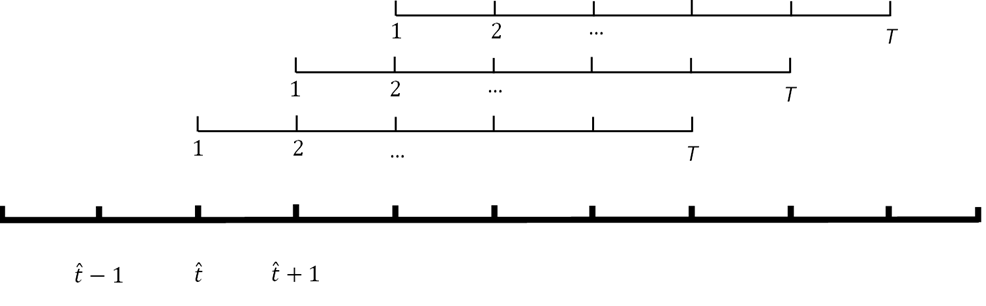}
\caption{Rolling-horizon framework} 
\label{fig:FiniteVSInfinite}
\end{center}
\end{figure}

To provide a decision policy for this infinite-horizon problem we propose a rolling-horizon approach, where at each time period we solve a finite-horizon problem that looks ahead $\Ti$ time periods and implement the first-period decisions obtained from this problem, as illustrated in Figure~\ref{fig:FiniteVSInfinite}. When the current period is $\tHat$, this finite-horizon model considers periods $\tHat, ..., \tHat +\Ti -1$. For notational convenience, when describing this model we shift all periods back by $\tHat-1$, so that the planning horizon is $\TI=\{1, ..., \Ti\}$.

In this section we formulate a finite-horizon multi-stage stochastic programming problem that we would ideally solve in each time period to make the current period decisions. This problem is a dynamic stochastic program with chance constraints to represent the service level requirements, and hence is intractable to solve exactly. In Section \ref{sec:approx} we discuss our proposed approximate solution strategies which are based on solving approximations of this finite-horizon problem.

\begin{figure}[ht]
\begin{center}
\includegraphics[scale=0.55]{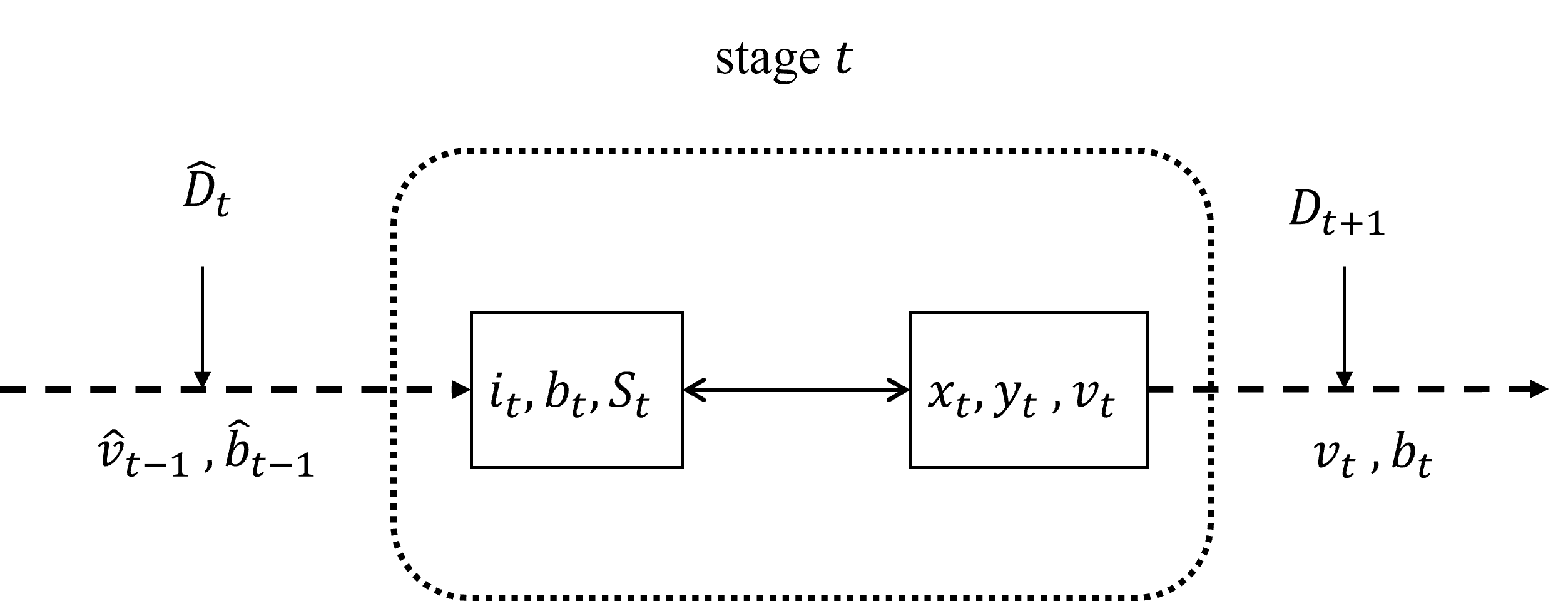}
\caption{Dynamics of decisions at each stage} 
\label{MultistageDynamics}
\end{center}
\end{figure}

Being at period $\ti=1$, given the state of the system, the model considers decisions for the $\Ti$ stages to guide the implementable first-stage ($\ti =1$) decisions that satisfy a joint service level in the next period, $\ti=2$.
Figure~\ref{MultistageDynamics} illustrates the dynamics of decisions at each stage $\ti$. 
First, the demand realizations $\hat{D}_{\ti \ka}$ for $\ka \in \KA$ are observed and we also know the initial state of the system, described by the on-hand inventory amounts $\hat{\Vi}_{\ti-1,\ka}$ and backlog amounts $\hat{\Bi}_{\ti-1,\ka}$ for $\ka \in \KA$. Based on this information, two sets of decisions are made.
The first set of decisions are the substitution decisions, which then imply the intermediate inventory and backlog amounts. 
The inventory of a product $\ka \in \KA$ can be used to satisfy demand of any product in the set $\Csub \subseteq \KA$, where we assume $\ka \in \Csub$ indicating that inventory of a product can certainly be used to meet its own demand. For each $\ka \in \KA$ we also define the set $ \Psub \subseteq \KA$ to be the set of products whose inventory can be used to meet demand of product $\ka$, and observe that for a pair of products $(\ka,\jey)$, $\jey \in \Csub$ if and only if $\ka \in \mathcal{K}^-_{\jey}$.
Thus, for each $\ka \in \KA$ and $\jey \in \Csub$, ${\Es}_{\ti \ka \jey}$ represents the amount of inventory of product $\ka$ that is used to meet demand or backlog of product $\jey$. Note that  $\Es_{\ti \ka \ka}$ corresponds to the amount of product $\ka$ which is used to satisfy its own demand. The substitution decisions, together with the demand, on-hand inventory, and backlog amounts then define the intermediate inventory and backlog amounts $\inv_{\ti \ka}$ and $\Bi_{\ti \ka}$ for $\ka \in \KA$, respectively. 
The second set of decisions are the setup and production decisions, which combined with the intermediate inventory  levels determine the \InvPos \ at the end of current period. For each product $\ka \in \KA$, $x_{\ti \ka}$ represents the production amount of product $\ka$, $y_{\ti \ka}$ is a binary variable indicating if a setup is done ($y_{\ti \ka} = 1$) or not ($y_{\ti \ka} = 0$), and $\Vi_{\ti \ka}$ represents the \InvPos\ at the end of this period (equivalently, the beginning of the next period). 
Note that the substitution and production decisions are made simultaneously, but our convention that demand is observed at the beginning of a period and the production lead time is one period implies that the production amounts decided in period $\ti$ can only be used to meet demand or fill backlog in the next period or later.  This is why we have two different inventory levels for each product, namely, $\inv_{\ti \ka}$ as the inventory level immediately after demand satisfaction, but before production, and $\Vi_{\ti \ka}$ as the \InvPos\ at the end of the period. 
The values of $\Vi_{\ti \ka}$ and $\Bi_{\ti \ka}$ for all $\ka \in \Ka$ are the inputs for the next period, describing the next state of the system.

For each $\ka \in \KA$ and $\jey \in \Csub$, $c^{\text{sub}}_{\ti \ka \jey }$ represents the cost incurred at period $\ti$ when a unit of product $\ka$ is used to meet unit of demand of product $\jey$. Typically, $c^{\text{sub}}_{\ti \ka \ka} = 0$ representing that there is no additional cost incurred when a product is used to meet its own demand.
 An inventory holding cost of $c^{\text{hold}}_{\ti \ka}$ is charged for each unit of product $\ka \in \KA$ held in inventory after the demand satisfaction in period $\ti$. The cost to produce a unit of product $\ka \in \KA$ in period $\ti$ is denoted by $c^{\text{prod}}_{\ti \ka}$. Furthermore, if a setup of product $\ka \in \KA$ is done in period $\ti$, a setup cost of $c^{\text{setup}}_{\ti \ka}$ is incurred.



\begin{table}[ht]
\centering
\caption{Notation for the mathematical model}
\begin{adjustbox}{width=1\textwidth,center=\textwidth}
\begin{tabular}{ll}
\toprule
{\textbf {Sets}} & {\textbf {Definition}} \\ \midrule
$\TI$  & Set of planning periods, indexed by $1, ... ,\Ti$ \\ 
$\KA$  & Set of products, indexed by $1, ... ,\Ka$ \\
$  \Csub$ 
& Set of products whose demand can be fulfilled by product $\ka$  \\
$ \Psub$ 
& Set of products that can fulfill the demand of product $\ka$  \\
\midrule 
{\textbf {Parameters}} & {\textbf {Definition}} \\ \midrule
$c^{\text{setup}}_{\ti \ka}$ & Setup cost for product $\ka$ in period $\ti$ \\ 
$c^{\text{hold}}_{\ti \ka}$  & Inventory holding cost for product $k$ in period $\ti$  \\ 
$c^{\text{sub}}_{\ti \ka \jey }$  & Substitution cost if product $\ka$  is used to fulfill the demand of product $\jey$  in period $\ti$  \\ 
$c^{\text{prod}}_{\ti \ka}$  & Production cost for product $\ka$ in period $\ti$  \\
$c^{\text{back}}_{\ti \ka}$  & Backlog cost for product $\ka$  in period $\ti$ \\
$\alpha$  & Minimum required joint service level \\ 
$M_{\ti \ka}$  & Maximum production of product $\ka$ in period $\ti$ \\ 
${D}_{\ti \ka}$ & Random variable representing demand for product $\ka $ in period $\ti$  \\ 
${D}^\text{Hist}_{\ti \ka}$ & Vector of random demands from period 1 to period $\ti$ for product $\ka $  \\ 
$\hat{\Vi}_{0,\ka} $&  The amount of initial inventory level for product $\ka$ \\
$\hat{\Bi}_{0,\ka} $&  The amount of initial backlog for product $\ka$   \\
$\mathbb{P} $&  The probability distribution of the demand process \\ 
\midrule
{\textbf {Decision variables}} & {\textbf {Definition}} \\ \midrule
$\y_{\ti \ka}$ & Binary variable which is equal to 1 if there is a setup for product $k$ at period $\ti$, 0 otherwise \\ 
$\x_{\ti \ka}$ & Amount of production for product $\ka$  at period $\ti$  \\ 
$\Es_{\ti \ka \jey}$ & Amount of product $\ka$  used to fulfill the demand of product $\jey$  at period $\ti$   \\
$\inv_{\ti \ka}$ & Amount of physical inventory for product $k$ immediately after the demand satisfaction for period $\ti$  \\
${\Bi}_{\ti \ka}$ & Amount of backlog for product $k$ at the end of period $\ti$  \\
${\Vi}_{\ti \ka}$ & The \InvPos\ for product $\ka$ at the end of period $\ti$ (beginning of period $\ti+1$)  \\
 \bottomrule
\end{tabular}
\end{adjustbox}
 \label{tab:Sub_parameters}
\end{table}

The product demands are modeled as a stochastic process, $D_{\ti \ka}$ for $\ti=1,\ldots,T$ and $\ka \in \Ka$, where $D_{\ti \ka}$ is a random variable representing the demand of product $\ka$ in period $\ti$. We use the notation $\hat{D}_{\ti \ka}$ to indicate a particular observed realization of this random variable.
$D^\text{Hist}_{\ti \ka}$ represents the random demand path from period 1 to period $\ti$ for product $\ka$, and $\hat{D}^\text{Hist}_{\ti\ka}$ denotes its realization (the history) until period $\ti$.

For notational convenience, when an index is dropped when referring to a parameter or decision variable, we are referring to the vector of all the parameters and decision variables over the range of that index. For instance,  $\hat{D}_{\ti} := (\hat{D}_{\ti 1}, \ldots, \hat{D}_{\ti \Ka})$ and likewise for ${\Es}_{\ti}$, $\inv_{\ti}$, ${\Bi}_{\ti}$, $\Vi_{\ti}$, etc.

We now present the finite-horizon chance-constrained multi-stage stochastic programming model.
Notation for different sets, parameters, and decision variables is summarized in Table~\ref{tab:Sub_parameters}. 
 We present a dynamic programming formulation where $F_t(\cdot)$ denotes the cost-to-go function at each period $\ti = 1, 2, ..., \Ti$ and is defined recursively as follows:
\begin{subequations}
\label{DynamicProgrammingNoden}
\begin{alignat}{2}
F_{t}(\hat{\Vi}_{\ti -1}, \hat{\Bi}_{\ti -1}, \hat{D}^\text{Hist}_\ti) = \min \ & \sum_{\ka \in \KA} \Bigl( c^{\text{setup}}_{\ti \ka}\y_{\ti \ka} + c^{\text{prod}}_{\ti \ka}\x_{\ti \ka} + c^{\text{hold}}_{\ti \ka} \inv_{\ti \ka}+  \sum_{\jey \in \Csub}c^{\text{sub}}_{\ti \ka \jey} \Es_{\ti  \ka \jey}  \Bigr) + && \quad  \notag \\*[0.1cm]
& \mathbb{E}_{D_{\ti+1}} \left[{F}_{\ti+1}(\Vi_{\ti}, \Bi_{\ti} , {D}^{\text{Hist}}_{\ti+1}) |{D}^{\text{Hist}}_{\ti} = \hat{D}^{\text{Hist}}_{\ti}\right] \label{eq:Dyn_g2_Sub_ST_Production_Flow} \\
\text{s.t.} \ & \x_{\ti \ka} \leq M_{\ti \ka} \y_{\ti \ka}  && \forall \ka  \in \KA \label{eq:Dyn_Sub_Setup} \\
& \sum_{\jey \in  \Psub} \Es_{\ti \jey \ka} + \Bi_{\ti  \ka}  = \hat{D}_{\ti \ka} + \hat{\Bi}_{\ti-1, \ka} && \forall \ka  \in \KA \label{eq:Bhat} \\
& \sum_{\jey \in  \Csub} \Es_{\ti \ka \jey} + \inv_{\ti \ka} = \hat{\Vi}_{\ti-1, \ka}  && \forall \ka  \in \KA  \label{eq:vhat} \\
& \Vi_{\ti \ka} = \inv_{\ti  \ka} + \x_{\ti  \ka}  && \forall \ka \in \KA \label{eq:vdef} \\
& \mathbb{P}_{D_{\ti+1}}\{ ({\Vi}_{\ti}, {\Bi}_{\ti} ) \in Q(D_{\ti+1} )| D^\text{Hist}_{\ti} = \hat{D}^\text{Hist}_{\ti} \} \geq \alpha \label{eq:SL} \\*[0.1cm]
& (\hat{\Vi}_{t-1},\hat{\Bi}_{t-1}) \in Q(\hat{D}_t) \Rightarrow \Bi_t = \mathbf{0}, && \label{eq:csm} \\
& {x}_{ \ti },  \Vi_{ \ti },  \inv_{ \ti } , {\Bi}_{ \ti } \in \mathbb{R}_+^{\Ka}, \Es_{\ti} \in \mathbb{R}_+^{\Ka \times \Ka}, \ {y}_{ \ti } \in \{0,1\}^{\Ka} \label{eq:Dyn_F_Sub_ST_bound1} 
\end{alignat} 
\end{subequations}
where $F_{\Ti+1}(\cdot) =0$ and $\mathbf{0}$ denotes the vector of zeros of appropriate dimension.

 The optimal value,  denoted by $F_t(\hat{\Vi}_{\ti -1}, \hat{\Bi}_{\ti -1}, \hat{D}^\text{Hist}_\ti)$, represents the optimal objective value from period $\ti$ to the end of the horizon given the initial \InvPos\ vector $\hat{\Vi}_{\ti-1}$, backlog vector $\hat{\Bi}_{\ti-1}$, and observed demand history $\hat{D}^{\text{Hist}}_\ti$. The objective \eqref{eq:Dyn_g2_Sub_ST_Production_Flow} is to minimize the current stage total cost (setup, production, holding, and substitution costs) plus the expected optimal costs from stages $\ti+1$ to the end of the horizon. 
Constraints (\ref{eq:Dyn_Sub_Setup}) are the setup constraints which ensure that if there is production of a product $\ka \in \KA$, the setup variable $y_{\ti \ka}$  takes the value 1. Here,  $M_{\ti \ka}$ is an upper bound on the maximum amount of product $\ka$ that would be produced in period $\ti$ in an optimal solution.
Constraints (\ref{eq:Bhat}) enforce that the current demand plus last period's backlog of each product is satisfied or it will be recorded as backlog $\Bi_{\ti \ka}$ for the next period.
Constraints (\ref{eq:vhat}) model the use of \InvPos\ $\hat{\Vi}_{\ti-1, \ka}$ of each product $\ka \in \KA$. It may be used to meet demand of any product in the set $\Csub$ or it will be recorded as intermediate inventory and combined with current period production $x_{\ti \ka}$ to yield the next period's \InvPos\ $\Vi_{\ti}$, as described in 
constraints (\ref{eq:vdef}). 

The {\it stock-out free set} $Q(D)$, defined for a vector of demands $D=(D_1,\ldots,D_{\Ka})$ plays an important role in constraints \eqref{eq:SL} and  \eqref{eq:csm}. This set represents the set of \InvPos\  and backlog vectors for which it is possible to avoid a stock-out of any product in the next period if the product demands are given by the vector $D$. Specifically, the set is defined as:
\begin{alignat}{2} Q(D) := \bigl\{ (\Vi,\Bi) \in \mathbb{R}_{+}^{\Ka} \times \mathbb{R}_+^{\Ka} :  & \ \exists \Es_{\ka \jey} \geq 0, \forall \ka \in \KA, \jey \in \Csub && \text{ such that } \notag \\
& \sum_{j \in \Psub} \Es_{\jey \ka } = D_{\ka} + \Bi_{\ka} \ && \forall \ka  \in \KA  \quad \text{and }\notag \\
 & \sum_{\jey \in  \Csub} \Es_{\ka \jey} \leq \Vi_{\ka} \ && \forall \ka  \in \KA \bigr\}  \label{eq:SetQ}
 \end{alignat}
 so that $(\Vi,\Bi) \in Q(D)$ if and only if it is possible to meet all product demands $D$ and backlogs  $\Bi$ using available inventory $\Vi$.
Thus, constraint~(\ref{eq:SL}) requires that there is sufficient inventory in the next period to avoid a stock-out with probability at least $\alpha$, where this probability is over the distribution of the next-stage's demand conditional on the current history $\hat{D}_t^{\text{Hist}}$. While this constraint ensures that there is always at least an $\alpha$ probability that stock-outs can be avoided in the next stage, the constraint by itself is not sufficient to enforce the $\alpha$ service level, due to the possibility to allow a stock-out to occur when making substitution decisions (e.g., to save costs) even though it might be feasible to avoid one. This is the purpose of constraint \eqref{eq:csm} -- it states that if a stock-out can be avoided in the current stage (i.e., $(\hat{\Vi}_{\ti-1},\hat{\Bi}_{\ti-1})$ is in the stock-out free set for the current demands), then a stock-out is not allowed (i.e., $\Bi_{\ti k} = 0$ for all $\ka \in \KA$). This constraint reflects a modeling assumption that the firm always wishes to avoid stock-outs when feasible, e.g., to avoid difficult-to-quantify costs such as loss of customer goodwill, an assumption we argue is consistent with the use of the $\alpha$ service level constraint. Although we do not pursue this possibility here, an alternative to constraint \eqref{eq:csm} would be to introduce decision variables that define a policy for determining when a stock-out will be allowed, and include the optimization of those decision variables as part of the formulation.

In order to assure that there is always a feasible solution to problem \eqref{DynamicProgrammingNoden} we assume that for each product $\ka \in \KA$ there is a product $\jey \in \Csub$ whose production limit $M_{\ti \jey}$ is large enough that it is always possible to produce enough in the current period to avoid stock-outs in the next period with the desired minimum probability.

\section{Approximate solution policies}
\label{sec:approx}

We now present our proposed approximate solution policies. 
As described in Section \ref{sec:form}, our approach is to solve a finite-horizon optimization model in each time period to make the current decisions given the current available inventory and backlog. 
Ideally, we would solve model \eqref{DynamicProgrammingNoden} and implement the solution from time period $1$.
Since this model is intractable due to its multi-stage nature and high-dimensional state space,
we instead propose to solve an approximation of this model and implement the decision that the approximation yields in the first period. We consider two types of approximations, one that is purely deterministic and one that incorporates a two-stage chance constraint to model the service level constraint. In both policies, the first step is to solve a model that determines if stock-outs can be avoided in the current stage (i.e., to enforce constraint \eqref{eq:csm}), which we describe in Section \ref{ss:csm}. The output of this model is then used to define constraints on backlog that are applied when we solve the finite-horizon approximate model to make decisions. These approximate models are described in Section \ref{ss:approx}. We describe a branch-and-cut algorithm for solving the two-stage chance-constrained model in Section \ref{ss:bnc}.

 \subsection{Stock-out determination}
 \label{ss:csm}

Given the current \InvPos\ vector $\hat{\Vi}_0$, backlog vector $\hat{\Bi}_0$, and observed demand $\hat{D}_1$, the first step is to determine whether a stock-out can be avoided in the first period, i.e., to test if $(\hat{\Vi}_0,\hat{\Bi}_0) \in Q(\hat{D}_1)$ as in \eqref{eq:csm}.
To this end, we solve the following linear programming (LP) model which minimizes the total backlog in the current period:  
\begin{subequations}
\label{Currentstage}
\begin{alignat}{2}
\min \ & \sum_{\ka \in \KA}  {\Bi}_{\ka}  \label{eq:Current_obj} \\
 \text{s.t.} \ &   \sum_{\jey \in  \Psub} \Es_{  \jey \ka} + \Bi_{\ka}  = \hat{D}_{1\ka} + \hat{\Bi}_{0 \ka} &&\qquad\forall \ka \in \KA       \label{eq:Current_inventory_tn}\\
&  \sum_{\jey \in  \Csub} \Es_{\ka \jey} \leq \hat{v}_{ 0 \ka} &&\qquad\forall \ka \in \KA      \label{eq:Current_Orderup}\\
& {\Bi} \in \mathbb{R}_+^{\Ka}, \Es \in \mathbb{R}_+^{\Ka \times \Ka}.
\label{eq:Current_bound2}
\end{alignat}
\end{subequations}
Constraints (\ref{eq:Current_inventory_tn}) guarantee that the demand and current period backlog is either satisfied or it will be backlogged in the next period. 
Constraints (\ref{eq:Current_Orderup}) limit the use of available inventory.
Given an optimal solution $(\Bi^*,\Es^*)$ of \eqref{Currentstage}, we set $\hat{\KA} = \{ \ka \in \KA: \Bi^*_{\ka}=0 \}$, which represents the set of products for which we are able to achieve zero backlog. In all our policies, we enforce $\Bi_{1 \ka} = 0$ for all $\ka \in \hat{\KA}$ so that we only allow backlog when it is impossible to avoid. In particular, if the optimal value of \eqref{Currentstage} is zero, then backlogging will not be allowed for any product, and hence the period will not experience a stock-out.
We stress that this model is only used to determine if we allow backlog at the end of the first period -- the actual substitution decisions are made after solving another model which we describe next.



\subsection{Approximate models} 
\label{ss:approx}
We now describe the deterministic approximation (Section \ref{ss:det}) and two-stage chance-constrained approximation (Section \ref{ss:cc}) that we propose to solve to make the current decisions in each time period.
The advantage of the deterministic approximation is that it is a MILP model, and hence is solvable by widely available MILP software. While the chance-constrained model is more difficult to solve, we find that it results in better policies, and can be solved efficiently with the branch-and-cut method we present in Section \ref{ss:bnc}.

\subsubsection{Deterministic approximation}
\label{ss:det}

The deterministic approximation is based on replacing all future uncertain demands with a deterministic estimate $\overline{D}_{\ka \ti}$ for $t=2,\ldots,T$ and $\ka \in \KA$, resulting in the following multi-period deterministic MILP model:
\begin{subequations}
\label{mod:Det}
\begin{alignat}{2}
\min \ & \sum_{\ti \in \TI} \sum_{\ka \in \KA} \Bigl( c^{\text{setup}}_{\ti \ka}y_{\ti \ka} + c^{\text{prod}}_{\ti \ka}x_{\ti \ka}+ c^{\text{hold}}_{\ti \ka}\inv_{\ti \ka} +\sum_{\jey \in  \Csub}c^{\text{sub}}_{\ti \ka \jey} \Es_{\ti \ka \jey}  \Bigr) 
& \label{eq:Sub_Det_obj} \\
\text{s.t.} \ & \  x_{\ti \ka} \leq M_{\ti \ka} y_{\ti \ka} &&  \forall \ti  \in \TI, \forall \ka \in \KA  \label{eq:Sub_Det_Setup}\\
  &  \sum_{\jey \in  \Psub} \Es_{1  \jey \ka} + \Bi_{1 \ka}  = \hat{D}_{1 \ka} + \hat{\Bi}_{0\ka} &&\forall \ka\in \KA       \label{eq:Det_inventory_1}\\
   &  \sum_{\jey \in  \Psub} \Es_{2 \jey \ka}  = \overline{D}_{\ka 2} + \Bi_{1 \ka} && \forall \ka\in \KA      \label{eq:Det_inventory_2} \\
   &  \sum_{\jey \in  \Psub} \Es_{\ti \jey \ka}  = \overline{D}_{\ka \ti}  &&\forall \ti \in \TI \setminus\{1,2\},\forall \ka\in \KA      \label{eq:Det_inventory_3} \\
&  \sum_{\jey \in  \Csub} \Es_{\ti \ka \jey} + \inv_{\ti \ka} = \Vi_{\ti-1 , \ka} &&\forall \ti  \in \TI
,\forall \ka \in \KA       \label{eq:Det_inventory_4}\\
& \Vi_{\ti \ka} = \inv_{\ti \ka} + \x_{\ti \ka}  \quad &&\forall \ti  \in \TI,\forall \ka \in \KA  \label{eq:Det_inventory_5} \\
&\Bi_{1 \ka} =0 &&\forall  \ka \in \hat{\KA} \label{eq:Backlog_determination}\\
&  {x}_{ \ti },  {\Vi}_{ \ti },  \inv_{ \ti} , {\Bi}_{ \ti } \in \mathbb{R}_+^{\Ka}, \Es_{\ti} 
\in \mathbb{R}_+^{\Ka \times \Ka}, {y}_{ \ti } \in \{0,1\}^{\Ka} &&\forall \ti \in \TI  \label{eq:Sub_FD_bound2}
\end{alignat}
\end{subequations}
The objective function (\ref{eq:Sub_Det_obj}) minimizes the total cost of setup, production, holding and substitution cost over the $T$ planning periods. Constraints (\ref{eq:Sub_Det_Setup}) guarantee that in each planning period, when there is positive production, there will be a setup. In case the production level of a product $k \in \KA$ is unconstrained in period $\ti \in \TI$, a sufficiently large value of $M_{\ti \ka}$ for use in constraint \eqref{eq:Sub_Det_Setup} can be computed as:
\begin{alignat}{2}
  &  M_{\ti \ka} =  \sum_{\jey \in  \Csub}\Bigl( \hat{\Bi}_{0 \jey}+ {\hat{D}}_{1 \jey} + \sum_{t=2}^T \overline{D}_{tj} \Bigr).  
  \label{eq:BigM_Deterministic_Appr}
  \end{alignat}
Constraints (\ref{eq:Det_inventory_1}) to (\ref{eq:Det_inventory_4}) are the inventory, backlog, and substitution balance constraints.
In constraints (\ref{eq:Det_inventory_4}) for $\ti=1$, $\Vi_{0 , \ka} := \hat{\Vi}_{0\ka}$.
Constraints~(\ref{eq:Det_inventory_1}) and \eqref{eq:Det_inventory_2} use the $b_{1k}$ variables, which according to (\ref{eq:Backlog_determination}) are only allowed to be positive when stock-out could not be avoided, as determined in the stock-out determination step. Constraints \eqref{eq:Det_inventory_2} and \eqref{eq:Det_inventory_3} do not use backlog variables for any $\ti > 1$, and hence this model enforces satisfaction of the deterministic estimates of demand in periods $\ti > 1$. 
Constraints (\ref{eq:Det_inventory_5}) define the available inventory after production. 

\begin{figure} [ht]
    \centering 
    \subfloat[\centering Average policy ]{{\includegraphics[width=5cm]{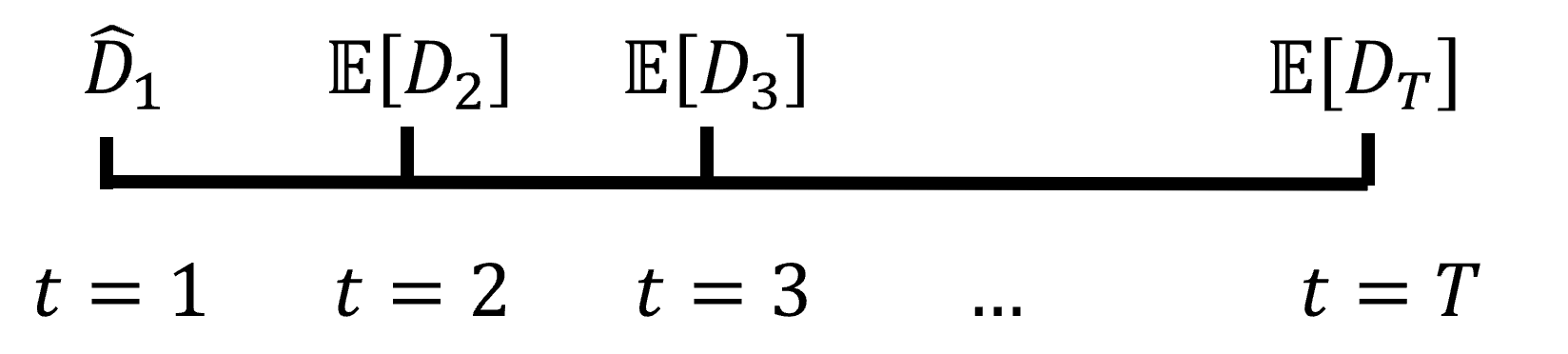} }}%
    \subfloat[\centering Quantile policy ]{{\includegraphics[width=5cm]{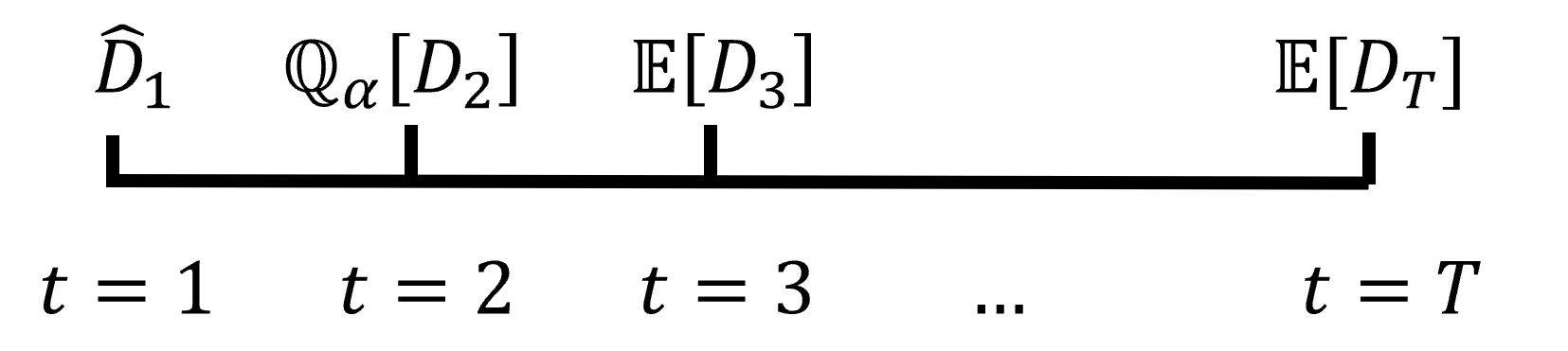} }}%
    \subfloat[\centering Chance-constraint policy ]{{\includegraphics[width=5.3cm]{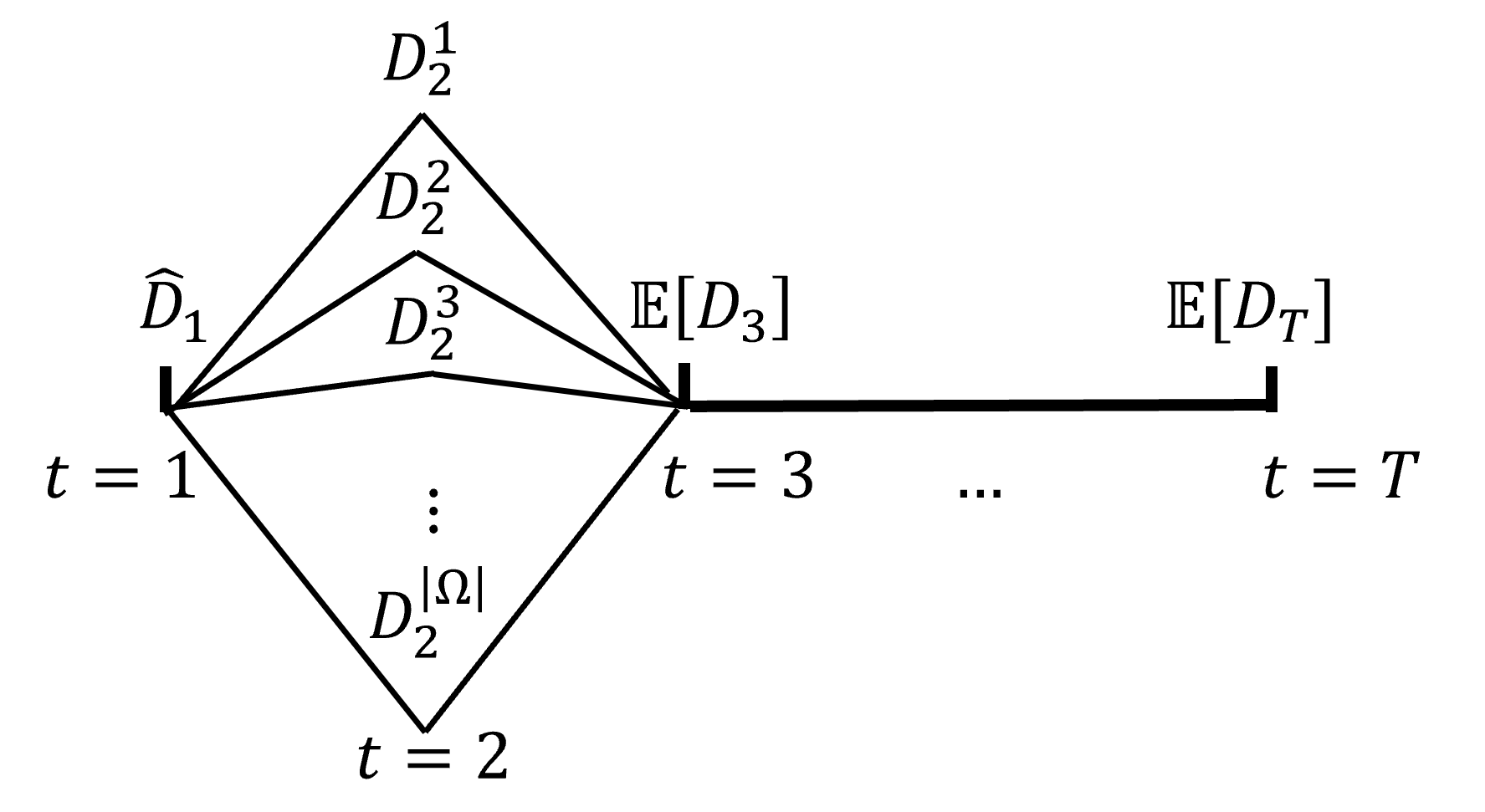} }}%
    \caption{Demand approximation in different decision policies}%
    \label{fig:DemandPolicy}%
\end{figure}

We consider two variations of the deterministic approximation based on using different estimates for the future demands.
In the first policy, which we refer to as the ``average policy", the expected value of demand is used for the deterministic approximation, specifically $\overline{D}_{\ti \ka} = \mathbb{E}[D_{\ti \ka}]$ for $\ti=2,\ldots,T$ and all $\ka \in \KA$. This average policy has little chance of meeting the service level constraints, since the decisions made in the current period are only planning for the expected demand of each product in the next period, so that the realized demand in the next period will often exceed the amount planned for.  We thus consider a second policy, which we refer to as the ``quantile policy", where the demand in the next immediate period is approximated by {\it the $\alpha$ quantile} of the future demand distribution for each product, and the demand in periods beyond that are approximated by their expected value. So, in this case $\overline{D}_{2 \ka} = \mathbb{Q}_{\alpha}[D_{2 \ka}] := \min\{q : \mathbb{P}(D_{2 \ka} \leq q) \geq \alpha\}$  and  $\overline{D}_{\ti \ka} = \mathbb{E}[D_{\ti \ka}]$ for $\ti=3,\ldots,T$ and all $\ka \in \KA$. Figure~\ref{fig:DemandPolicy} illustrates the demand pattern for the average and quantile policies in sub-figures (a) and (b), respectively.

\subsubsection{Two-stage chance-constrained approximation}
\label{ss:cc}

The deterministic policies that we explained in the previous section do not consider the service level constraint explicitly. To ensure the service level is met, we introduce a joint chance constraint that requires current ordering decisions be sufficient to ensure that all demand can be met in the next period in at least $\alpha$ fraction of the possible demand outcomes in the next period. To model the chance constraint, we approximate the joint distribution of product demands in period $\ti=2$ using a finite set of equally likely joint demand scenarios $D^{\m}_{2}$ for $\m \in \EM$, where $D^{\m}_{2\ka}$ denotes the demand of product $\ka \in \KA$ in period $2$ under scenario $\m \in \EM$. The finite set of scenarios could be obtained, for example, via Monte Carlo sampling \citep{luedtke2008sample}. 

Given a discrete approximation of the next period's demand distribution, it is possible to formulate the service level constraint \eqref{eq:SL} explicitly by introducing additional variables to represent which scenarios will not have a stock-out in the next stage.
However, just modeling this constraint is likely to lead to poor policy performance because it would ignore the cost of the next-stage substitution that is implicitly planned for when enforcing the chance constraint.  Specifically, having only the chance constraint would ensure that with high probability there exists substitutions in the next stage that
can meet all demands, but would ignore the cost of those substitutions. 
This would in turn lead to decisions that require potentially costly product substitutions to avoid stock-outs. To address this issue, we propose a model that {\it both} enforces a service level constraint {\it and} approximates the cost of the decisions in period 2 for each single scenario $\m \in \EM$. To preserve compactness of the model, we continue to approximate the demand in periods $\ti \geq 3$ with the expected values $\mathbb{E}[D_{\ti\ka}]$ for $\ka \in \KA$. The demand pattern used in this policy in depicted in sub-figure (c) in Figure~\ref{fig:DemandPolicy}.

This model uses new decision variables representing the decisions that will be made in each scenario in period $2$: $I'^{\m}_{2\ka }$ and $B'^{\m}_{2 \ka}$ denote the inventory and backlog at period 2 for product $\ka \in \KA$ under scenario $\m \in \EM$, respectively. 
The variable $\Es'^{\m}_{2 \jey \ka}$ represents the amount of product $\jey$ used to meet demand of product $\ka$ under scenario $\m \in \EM$ in period 2. 
The proposed model is:
\begin{subequations}
\label{mod:ChanceConstraintpolicy}
\begin{alignat}{2}
\min \ &
\sum_{\ka \in \KA} \Bigl( c^{\text{setup}}_{ 1 \ka}y_{ 1 \ka} + c^{\text{prod}}_{ 1 \ka}x_{ 1 \ka}+ \sum_{\jey \in  \Csub}c^{\text{sub}}_{1 \ka \jey} \Es_{1 \ka \jey}  + c^{\text{hold}}_{ 1 \ka}\inv_{ 1 \ka} \Bigr) 
+ \notag \\
&\sum_{\ka \in \KA} \Bigl( c^{\text{setup}}_{ 2 \ka}y_{2  \ka} + c^{\text{prod}}_{ 2 \ka}x_{2  \ka}+  c^{\text{hold}}_{2  \ka}I_{2  \ka} + c^{\text{back}}_{2 \ka} \Bi_{2  \ka} + \frac{1}{\Em} \sum_{\m \in \EM} \sum_{\jey \in  \Csub}  c^{\text{sub}}_{2\ka\jey}   \Es'^{\m}_{2\ka\jey } \Bigr) && + \notag \\
& \sum_{t =3 }^{T} \sum_{\ka \in \KA} \Bigl( c^{\text{setup}}_{\ti \ka}y_{\ti \ka} + c^{\text{prod}}_{\ti \ka}x_{\ti \ka}+ \sum_{\jey \in  \Csub}c^{\text{sub}}_{\ti \ka \jey} \Es_{\ti \ka \jey}  + c^{\text{hold}}_{\ti \ka}\inv_{\ti \ka} \Bigr) 
 \label{eq:Sub_Roll_obj_ext} 
 \end{alignat}
 \begin{alignat}{2}
\text{s.t.} \  & \ x_{\ti \ka} \leq M_{\ti \ka} y_{\ti \ka} &&  \forall \ti  \in \TI   ,\forall \ka \in \KA  \label{eq:Sub_cc_Setup}\\
  &  \sum_{\jey \in  \Psub} \Es_{1  \jey \ka} + \Bi_{1 \ka}  = \hat{D}_{1 \ka} + \hat{\Bi}_{0 \ka} &&\forall \ka \in \KA       \label{eq:Det_backorder_tn}\\
 &  \sum_{\jey \in  \Psub} \Es'^{\m}_{2 \jey \ka} + \Bi'^{\m}_{2  \ka }  = D^ {\m}_{\ka} + \Bi_{1  \ka} &&\forall \ka \in \KA, \forall \m \in \EM     \label{eq:Det_backorder_tnp} \\
    &  \sum_{\jey \in  \Psub} \Es_{\ti \jey \ka}   = \mathbb{E}[{D}_{\ti \ka}] + {\Bi}_{\ti -1 ,\ka} &&\forall \ti  \in \TI, \ti \geq 3,\forall \ka \in \KA     \label{eq:Det_backorder_ext} \\
   &  \sum_{\jey \in  \Csub} \Es_{1 \ka \jey} + \inv_{ 1 \ka} = \hat{v}_{0  \ka} &&\forall \ka \in \KA       \label{eq:Det_inventory_tn}\\
&  \sum_{\jey \in  \Csub} {\Es}'^{\m}_{2 \ka \jey} + \inv'^{\m}_{2  \ka} = \Vi_{1  \ka} &&\forall \ka \in \KA, \forall \m \in \EM       \label{eq:Det_inventory_tnp}\\   
&  \sum_{\jey \in  \Csub} \Es_{\ti \ka \jey} + \inv_{\ti \ka} = v_{\ti-1 , \ka} &&\forall \ti  \in \TI, \ti \geq 3,\forall \ka \in \KA       \label{eq:Det_inventory_ext}\\
& \Vi_{\ti \ka} = \inv_{\ti \ka} + x_{\ti \ka}  \quad &&\forall \ti  \in \TI,\forall \ka \in \KA  \label{eq:Chance_OrderUpTo} \\
&\Bi_{1 \ka} =0 &&\forall  \ka \in \hat{\KA} \label{eq:Backlog_determination_CC}\\
& \frac{1}{\Em} \sum_{\m \in \EM} \inv'^{\m}_{2 \ka } = \inv_{2 \ka} &&\forall \ka \in \KA  \label{eq:Average_Inventory} \\
& \frac{1}{\Em} \sum_{\m \in \EM} \Bi'^{\m}_{2 \ka} = \Bi_{2 \ka} &&\forall \ka \in \KA  \label{eq:Average_Backlog}\\
& \sum_ {\m \in \EM}  \mathds{1} {\{ ({\Vi}_{1}, {\Bi}_{1} ) \in Q(D^ {\m}_{2} ) \}} \geq \lceil \alpha |\EM|  \rceil \label{eq:chance_ServiceLevel}\\
&{x}_{ \ti },  {\Vi}_{ \ti },  {\inv}_{ \ti} , {\Bi}_{ \ti } \in \mathbb{R}_+^{\Ka}, {\Es}_{\ti} 
\in \mathbb{R}_+^{\Ka \times \Ka}, {y}_{ \ti } \in \{0,1\}^{\Ka} \qquad &&\forall \ti  \in \TI  \label{eq:Sub_FD_bound3} \\
& {\inv}'^{\m}_{ 2} , {\Bi}'^{\m}_{ 2 }, {\Es}'^{\m}_{2} \in \mathbb{R}_+^{\Ka}, && \forall \m \in \EM
\end{alignat}
  \end{subequations}
The objective function in (\ref{eq:Sub_Roll_obj_ext}) is broken into three parts representing the cost in period 1, the cost of period $2$, and the cost of periods $3$ to $\Ti$. The cost is the same as the deterministic approximation in periods 1 and $t \geq 3$. In period $2$, the substitution cost is defined for each scenario separately, and the average substitution cost over all scenarios is included in this period's cost. Another key difference of the cost in period 2 is the presence of a backlog penalty term, with backlog ``cost'' parameter $c^{\text{back}}_{2 \ka}$ for $\ka \in \KA$. This term is included to encourage the decisions the model selects for the different scenarios in period 2 (the $\Es'^{\m}_{2},\inv'^{\m}_2,\Bi'^{\m}_2$ variables) to match the decisions that would actually be made when that period occurs and the stock-out determination step is applied to enforce that there are no stock-outs unless they cannot be avoided. To encourage this match, we suggest to select the backlog ``cost'' parameters $c^{\text{back}}_{2 \ka}$ so that the fraction of scenarios in which the backlog variables $\Bi'^m_2$ are zero across all products roughly approximates the desired service level. We stress that the values $c^{\text{back}}_{2 \ka}$ should be considered as a parameter of the policy, and are not meant to reflect actual backlog costs.

Constraints~(\ref{eq:Sub_cc_Setup}) are production setup constraints. In the case when production for a product $\ka$ does not have a given capacity, the $M_{\ti \ka}$ values can be set as
\begin{alignat*}{2}
  &  M_{\ti \ka} =  \sum_{\jey \in  \Csub}\left( \hat{\Bi}_{0 \jey}+  {\hat{D}}_{1 \jey} + \max_{\m \in \EM} {D}^{\m}_{2 \jey}+ \sum_{\ti =3}^\Ti \mathbb{E}[{D}_{\ti \jey}]\right)  &&\qquad \forall \ti \in \TI ,\forall \ka \in \KA .   
  \end{alignat*}
Constraints 
\eqref{eq:Det_backorder_tn}-\eqref{eq:Det_backorder_ext}
assure that demand plus carried over backlog are met or backlog is recorded in periods 1,2, and $t\geq3$, respectively. In period 1, the current observed demands and backlogs are what must be met (i.e., are in the right-hand side). In period 2, the demands for each different scenario are used. In periods $t\geq 3$, we use the expected demands.
Constraints \eqref{eq:Det_inventory_tn}-\eqref{eq:Det_inventory_ext} relate the available inventory in each period with how it is used and the inventory carried over to the next period.
Constraints (\ref{eq:Det_inventory_tn}) use the current available inventory $\hat{v}_{0\ka}$ for the current period constraint and restrict the substitution decisions and ending inventory accordingly. Constraints \eqref{eq:Det_inventory_tnp} consider analogous constraints in period 2, but do so for each scenario $\m \in \EM$, whereas constraints (\ref{eq:Det_inventory_ext}) present the analogous constraints for periods $t \geq 3$.
Constraints (\ref{eq:Average_Inventory}) and (\ref{eq:Average_Backlog}) define the variables $\Bi_{2\ka}$ and $\inv_{2 \ka}$ to be the averages of the $\Bi'^{\m}_{2 \ka}$ and $\inv'^{\m}_{2 \ka}$ variables over the set of scenarios $\m \in \EM$, respectively. These averaged variables are used in the inventory and backlog balance constraints for period $3$, and hence these constraints provide a critical link between the scenario variables used in period 2 to model the costs in different scenarios.
Finally, constraint (\ref{eq:chance_ServiceLevel}) represents the service level constraint. In this constraint, $\mathds{1}(\cdot)$ is an indicator that takes the value $1$ when the argument is true, and 0 otherwise. Thus, this constraint enforces that  $(\Vi_1,\Bi_1) \in Q(D_{2}^\m)$ (and hence $(\Vi_1,\Bi_1)$ is sufficient to meet demands $D_{2}^{\m}$ in period 2) in at least $\alpha$ fraction of the scenarios $\m \in \EM$. 
The constraint (\ref{eq:chance_ServiceLevel})
is not written in a form that can be given directly to a solver. In the next section we describe a branch-and-cut algorithm that can be used to solve the model with this constraint enforced.

  \subsection{Solving the chance-constrained model}
  \label{ss:bnc}
  
  We now discuss how to solve the proposed model \eqref{mod:ChanceConstraintpolicy}.
To do so, we define the binary variables $\Zed_\m$ for $\m \in \EM$ to model the indicator functions in  \eqref{eq:chance_ServiceLevel} and replace \eqref{eq:chance_ServiceLevel} with
\begin{equation}
\label{eq:Child_Service}
\sum_{\m \in \EM} \Zed_\m \leq \lfloor (1-\alpha) |\EM| \rfloor. 
\end{equation}
We must then enforce
\begin{equation}
\label{eq:zlogic}
\Zed_\m=0 \ \Rightarrow \ (\Vi_1,\Bi_1) \in Q(D_{2}^\m) \quad \forall \m \in \EM
\end{equation}
so that, for each scenario $\m \in \EM$, if $\Zed_\m = 0$ then the ending  available inventory $\Vi_1$ and backlog $\Bi_1$ are adequate to meet demands in period 2 without backlogging.
We present two options for enforcing the constraints \eqref{eq:zlogic}. 

\subsubsection{Extensive form}
\label{sec:ext}

In the first approach for enforcing the constraints \eqref{eq:zlogic}, we introduce variables
$\Bc^{\m}$ to represent the vector of backlog decisions and $\Sc^{\m}$ to represent the vector of substitution decisions in period 2 in scenario $\m$. The logical constraint \eqref{eq:zlogic} is then enforced with the following constraints:
  \begin{subequations}
\label{mod:chanceconstraintBigM}
  \begin{alignat}{2}
  & \Bc^ {\m}_{\ka } +\sum_{\jey \in \Psub} \Sc^{\m}_{\jey \ka} = D^ {\m}_{\ka}  +\Bi_{1 \ka}  && \qquad \forall \m \in \EM, \forall \ka  \in \KA \label{eq:Backorder_Child}\\
  & \sum_{\jey \in \Csub} \Sc^{\m}_{\ka \jey} \leq \Vi_{1 \ka}   &&\qquad \forall \m \in \EM, \forall \ka  \in \KA \label{eq:OrderUptoLevel_Child}\\
&  \Bc^ {\m}_{\ka} \leq \overline{M}^{\m}_{\ka }\Zed_{ \m}  && \qquad \forall \m \in \EM , \forall \ka  \in \KA     \label{eq:Child_Service_1}
 \end{alignat}
\end{subequations}
Constraints~\eqref{eq:Backorder_Child} and \eqref{eq:OrderUptoLevel_Child} define the backlog and substitution for each scenario $\m$. Constraints~(\ref{eq:Child_Service_1}) guarantee that if $\Zed_{ \m} =0$ then the backlog variables $\Bc^{\m}_{\ka}=0$ for all $\ka \in \KA$, so that the other constraints then enforce $(\Vi_1,\Bi_1) \in Q(D^{\m}_2)$. In (\ref{eq:Child_Service_1}), the $\overline{M}$ values are defined as:
\begin{alignat}{2}
 & \overline{M}^{\m}_{\ka }=  D^{\m}_{\ka} + \hat{D}_{1 \ka} + \hat{\Bi}_{0 \ka}  && \qquad \qquad \forall \m \in \EM , \forall \ka  \in \KA.      \label{eq:BigM_Child} 
 \end{alignat}
The variables $\Bc^{\m}$ and $\Sc^{\m}$ serve a similar role as the variables $\Bi'^{\m}_{2}$ and $\Es'^{\m}_{2}$  as described in Section \ref{ss:cc} in that they also represent backlog and substitution decisions in period 2. The difference is that the variables $\Bc^{\m}$ and $\Sc^{\m}$ are used to model the service level constraint, whereas the variables $\Bi'^{\m}_{2}$ and $\Es'^{\m}_{2}$ are used to approximate the cost of the decisions in period 2. Our next approach, which we find is computationally much more efficient than using \eqref{mod:chanceconstraintBigM}, does not introduce the variables $\Bc^{\m}$ and $\Sc^{\m}$.

\newcommand{\vsol}{\hat{\Vi}}
\newcommand{\zsol}{\hat{\Zed}}
\newcommand{\bsol}{\hat{\Bi}}
\newcommand{\pisol}{\pi}
\newcommand{\betasol}{\beta}
\newcommand{\optpi}{\hat{\pi}}
\newcommand{\optbeta}{\hat{\beta}}

\subsubsection{Branch-and-cut algorithm}
\label{sec:bnc}

The second approach for enforcing 
\eqref{eq:zlogic}
is to tailor the branch-and-cut algorithm proposed in \citep{luedtke2014branch} to this problem. 
In this approach, a master problem that includes the $\Zed_\m$ variables and the cardinality constraint 
\eqref{eq:Child_Service} (but not the $\Bc^{\m}$ and $\Sc^{\m}$ variables or constraints \eqref{mod:chanceconstraintBigM}) is constructed and cuts are iteratively added to it to enforce the logical constraints \eqref{eq:zlogic}. 

Assume we have solved a master problem and obtained a solution with $(\zsol,\vsol_1,\bsol_1)$ as the values for $(\Zed,\Vi_1,\Bi_1)$. Note that this solution may or may not satisfy the integrality constraints (e.g., if we have solved an LP relaxation of the master problem). Given a demand scenario $\m \in \EM$ with $\zsol^\m < 1$, our task is to assess if $(\vsol_1,\bsol_1) \in Q(D^\m_2)$, and if not, attempt to generate a cut to remove this solution. In the case of an integer feasible solution, we will always be able to do so when $(\vsol_1,\bsol_1) \notin Q(D_2^\m)$.

We can test if given $(\vsol_1,\bsol_1) \in Q(D_2^\m)$ by solving the following LP:
\begin{subequations}
\label{eq:feaschecklp}
\begin{alignat}{3}
V_\m(\vsol_1,\bsol_1) :=  \min_{w,\Sc^\m} \ & \sum_{\ka  \in \KA} w_\ka \\
    \text{s.t. } & \sum_{\jey \in  \Psub} \Sc^\m_{\jey \ka } + w_\ka = D^{\m}_{2\ka} + \bsol_{1\ka} \quad  &&\forall \ka  \in \KA && \label{eq:dual1}  \\
    &\sum_{\jey \in  \Csub} \Sc^\m_{\ka \jey} \leq \vsol_{1\ka} \ && \forall \ka  \in \KA&& \label{eq:dual2} \\
    & \ w \in \mathbb{R}_+^{\Ka}, \Sc^\m \in \mathbb{R}_+^{\Ka \times \Ka} 
\end{alignat}
\end{subequations}
By construction, $(\vsol_1,\bsol_1) \in Q(D_2^\m)$ if and only if $V_\m(\vsol_1,\bsol_1) \leq 0$, which means that there is no backlog for this scenario. Let $\pi$ and $\beta$ be the vectors of dual decision variables associated with constraints \eqref{eq:dual1} and \eqref{eq:dual2}, respectively, and let $\Pi$ be the set of dual feasible solutions. Observe that $\Pi$ is independent of $\m$ and $(\vsol_1,\bsol_1)$. 
Thus, for any $(\pisol,\betasol) \in \Pi$ and for any $\m \in \EM$, weak duality implies that the cut
\begin{equation}
\label{eq:basecutv0}
 \sum_{\ka  \in \KA} \pisol_k (D^{\m}_{2\ka} + \Bi_{1\ka}) + \sum_{\ka  \in \KA} \betasol_\ka \Vi_{1\ka} \leq 0
\end{equation}
is a valid inequality for $Q(D^\m_2)$. 
Since this inequality must hold whenever $\Zed_{\m}=0$, the inequality
\begin{equation}
\label{eq:basecut}
\sum_{\ka  \in \KA} \pisol_k \Bi_{1\ka} + \sum_{\ka  \in \KA} \betasol_\ka \Vi_{1\ka} \leq -  \sum_{\ka  \in \KA} \pisol_k D^{\m}_{2\ka} +  \bar{M}^{\betasol,\pisol}_{\m} \Zed_{\m} 
\end{equation}
is valid for suitably chosen (large enough) constant $\bar{M}_{\m}^{\betasol,\pisol}$.
Furthermore, if $(\pisol,\betasol)$ is taken to be the optimal dual solution to \eqref{eq:feaschecklp} for a given $(\vsol_1,\bsol_1)$ and scenario $\m \in \EM$, then if
$\zsol_{\m}=0$ and $V_\m(\vsol_1,\bsol_1) > 0$ then the corresponding cut is violated by $(\zsol, \vsol_1,\bsol_1)$, and hence is sufficient for  cutting off this solution whenever it violates \eqref{eq:zlogic}.

 We next discuss how to choose $\bar{M}_{\m}^{\betasol,\pisol}$ in \eqref{eq:basecut} and use this inequality to derive a family of strong valid inequalities that can be used to improve the LP relaxation. First, for each 
 $\m \in \EM$, define
\[ h_{\m}(\pisol,\betasol) = \sum_{\ka  \in \KA} \pisol_\ka D^{\m}_{2\ka}. \]
Using this notation in \eqref{eq:basecut}, 
we conclude that
\[ \sum_{\ka  \in \KA} \pisol_k \Bi_{1\ka} + \sum_{\ka  \in \KA} \betasol_\ka \Vi_{1\ka} \leq -  h_{\m}(\pisol,\betasol) \]
must be satisfied whenever $z_{\m} = 0$. We then sort the values $\{ h_{\m}(\pisol,\betasol) : \m \in \EM\}$ to obtain a permutation $\sigma$ of $\EM$ which satisfies:
\[ h_{\sigma_1}(\pisol,\betasol) \geq h_{\sigma_2}(\pisol,\betasol)  \geq \cdots \geq h_{\sigma_{|\EM|}}(\pisol,\betasol) .
 \]
Then, letting $p = \lfloor (1-\alpha) |\EM| \rfloor$, the followings 
are valid for the master problem~\citep{luedtke2014branch}:
\[ \sum_{\ka  \in \KA} \betasol_k \Vi_{1k} + \sum_{\ka  \in \KA} \pisol_k \Bi_{1k} \leq
  - h_{\sigma_i}(\pisol,\betasol) + \left(h_{\sigma_i}(\pisol,\betasol) - h_{\sigma_{p+1}}(\pisol,\betasol)\right)\Zed_{\sigma_i}, \quad \forall i=1,\ldots, p \]
and hence the coefficient on $\Zed_{\sigma_i}$ represents a valid value of $\bar{M}_{\omega}^{\betasol,\pisol}$ for $\omega = \sigma_i$ and $i=1,\ldots,p$.

A family of additional valid inequalities can be obtained by then applying {\it mixing inequalities} \citep{gunluk2001mixing,luedtke2014branch}. Given a subset $T = \{t_1,t_2,\ldots,t_{\ell}\} \subseteq \{\sigma_1,\sigma_2,\ldots,\sigma_p\}$ with $t_1 < t_2 < \cdots < t_{\ell}$ and defining $h_{t_{\ell+1}}:= h_{\sigma_{p+1}}$, the inequality
\begin{equation}
\label{eq:mixed}
  \sum_{\ka  \in \KA} \betasol_k \Vi_{1k} + \sum_{\ka  \in \KA} \pisol_k \Bi_{1k} \leq 
  - h_{t_1}(\pisol,\betasol) + \sum_{i=1}^{\ell} \left(h_{t_i}(\pisol,\betasol) - h_{t_{i+1}}(\pisol,\betasol)\right)\Zed_{t_i}
\end{equation}
is valid for the master problem. Although the number of such inequalities grows exponentially with $p$, there is an efficient algorithm for finding a most violated inequality~\citep{gunluk2001mixing} for given $(\zsol,\bsol_1,\vsol_1)$, which we describe for completeness in Algorithm~\ref{alg:mostviolated} in Appendix \ref{APP1}. 

Thus, given a solution $(\zsol,\bsol_1,\vsol_1)$ of the master problem, we proceed as follows to search for a cut. For any scenario $\m \in \EM$ with $\zsol_{\m} < 1$, we solve problem \eqref{eq:feaschecklp} to obtain a dual solution $(\optpi,\optbeta)$. If $V_\m(\vsol_1,\bsol_1) > 0$, we compute $h_{\omega'}(\optpi,\optbeta)$ for all $\m' \in \EM$ and finally search for a most violated inequality of the form \eqref{eq:mixed} and add it to the master problem, if violated. Within the branch-and-bound algorithm for solving the master problem, at the root node (i.e., the initial relaxation before branching) we carry out this process for any $\m \in \EM$ with $\zsol_{\m} < 1$ in the LP relaxation. At other solutions obtained in the branch-and-bound search, we only attempt to generate cuts when the solution $\zsol$ is integer-valued (and hence only for scenarios $\m$ with $\zsol_{\m}=0$). This is sufficient to guarantee that only solutions that satisfy \eqref{eq:zlogic} are accepted as feasible within the search process, thus leads to a correct solution. We refer to \citep{luedtke2014branch} for more details of the convergence analysis for this algorithm.

\section{Computational experiments}
\label{sec:comp}
We next report the results of our computational study which illustrate the ability of the proposed method to solve the chance-constrained model, demonstrate the benefits of the chance-constrained model-driven policy over policies based on deterministic approximations, and explore the benefits of substitution.

\subsection{Instance generation}

We generate a variety of test instances using \citep{rao2004multi} and \citep{hsu2005dynamic} as guidance for choosing substitution related parameters, and \citep{helber2013dynamic} for choosing the lot-sizing related parameters.
Table \ref{tab:BaseSensitivity} presents the key parameters we use to define an instance, their base value and the range of values we consider for this parameter when creating instances with different characteristics. We use $K=10$ products in all our tests. For the service level target $\alpha$, we range this between 80\% and 99\%, with 95\% as the base case.  The parameters $\eta$, $\tau$, $\rho$, and $TBO$  are used to calculate the cost parameters as described in Table
~\ref{tab:Sub_FD_parameters}. Parameter $\eta$ affects the relative difference in cost between the different products. Parameter $\tau$ impacts the cost of substitution (higher $\tau$ means substitution is more costly). The parameter $\rho$ determines the holding cost relative to the production cost, and the parameter TBO (time between orders) controls the relative setup cost. 

\begin{table}[ht] 
\small
\begin{minipage}[c]{0.65\linewidth} 
\begin{tabular}{lll}
\toprule
{Parameter} & Base Case & Variation \\ \midrule
$\Ka$   & 10 & \\ 
$\eta$  &   0.2 & 0.1, 0.2, 0.5   \\ 
$\tau$  &   1.5 & 1, 1.25, 1.5, 1.75, 2, 2.5   \\ 
$\rho $  &  0.05 & 0.02, 0.05, 0.1, 0.2 ,0.5  \\ 
$ TBO $  &   1 & 1, 1.25, 1.5, 1.75, 2   \\ 
$ \alpha $  &95\% & 80\%, 90\%, 95\%, 99\%  \\ 
\bottomrule 
\end{tabular} 
\caption{Base case and sensitivity analysis parameters} 
\label{tab:BaseSensitivity}
\end{minipage}%
\begin{minipage}[c]{0.35\linewidth} 
\begin{tabular}{ll}
\toprule
$\bar{c}^{\text{prod}}_{\ti \ka}$  & $1+\eta \times(\Ka-\ka)$ 
\\
$c^{\text{sub}}_{\ti \ka \jey }$  & $\max\{ 0,\tau \times (\bar{c}^{\text{prod}}_{ \ti \ka} - \bar{c}^{\text{prod}}_{\ti \jey}) \}$ 
\\ 
$c^{\text{hold}}_{\ti \ka}$  & $\rho \times \bar{c}^{\text{prod}}_{\ti \ka} $ 
\\ 
$c^{\text{setup}}_{\ti \ka}$ & $\mathbb{E}[{D_{\ti \ka}}]
\times TBO^2 \times c^{\text{hold}}_{\ti \ka} /2$ \\ 
\bottomrule
\end{tabular}
\caption{Cost parameters} 
\label{tab:Sub_FD_parameters} 
\end{minipage} 
\end{table}

In terms of the allowed substitution between products, we assume the products are ordered such that product 1 is the highest quality and product 10 is the lowest quality. In our base case, we assume that product $k$ can be used to meet demand of product $j$ if $k \leq j$ (so it is a higher quality product) and $k \geq j-3$ (so it is not more than three levels higher in the ranking).
Observe that when $\tau=1$ the cost of substituting a unit of product $k$ for a unit of product $j$ is exactly equal to the difference in production costs for these products. Thus, $\tau=1$ is a natural minimum value for this parameter in order to reflect the difference in production costs when substitution is performed, whereas larger values of $\tau$ represent the desire of a firm to limit the use of substitution, e.g., for business policy reasons. 

In every model that we solve, we enforce that the end-of-horizon backlog is zero for all products and hence the total amount of production is the same regardless of the model used. Since differences in production costs attributed to substitution are recorded in the substitution cost as described in the previous paragraph, we set all production costs $c^{\text{prod}}_{\ti \ka}$ equal to zero for all $\ka \in \KA$ and $\ti \in \TI$  in our experiments in order to exclude this cost from the cost comparison since it is constant across all policies. The parameters $\bar{c}^{\text{prod}}_{\ti \ka}$ in Table
 \ref{tab:Sub_FD_parameters} are used only to determine the values of the parameters $c^{\text{sub}}_{\ti \ka \jey }$ and $c^{\text{hold}}_{\ti \ka}$  as described in the table.


Recall that the chance-constrained model uses an artificial backlog cost on the backlog variables in period $2$ which needs to be tuned. We found that setting this parameter equal to the maximum possible cost of substitution yields reasonable results. Thus, we set
\begin{equation}
c^\text{back}_{ 2 \ka} = \max_{l \in \KA ,\jey \in \Psub }c^{\text{sub}}_{2 \jey l  }  \qquad \forall \ka  \in \KA.       \label{eq:backCost}
\end{equation}

The demands are assumed to be independent across different products, but demands for each product follow an auto-regressive (AR) model ~\citep{jiang2017production} which induces correlation in demand across time periods: 
\begin{equation}
D_{ \ti+1,\ka} = C + AR_1 D_{\ti \ka} + AR_2 \epsilon_{\ti+1,\ka }  \qquad \forall \ka  \in \KA , \forall \ti \in \TI,       \label{eq:AR1}
\end{equation}
where $C, AR_1,$ and $AR_2$ are parameters of the model, and $\epsilon_{ \ti+1,\ka}$ is a random noise with normal distribution with the mean of 0 and standard deviation of 1. 
In our data sets, we use $C = 20$, $AR_1 = 0.8$, and $AR_2 = 10$. Note that the expected demand for each product in each period is equal to $C/(1-AR_1)=100$. 

As we have no production in the first period, we assume that the demand in the first period is zero, otherwise, if there is no initial inventory, the service level constraint will not be satisfied. We initialize the AR data generation procedure with $D_{0 \ka}=C/(1-AR_1)$ (the expected demand), and then use \eqref{eq:AR1} with randomly generated values of $\epsilon_{\ti+1,\ka}$ to determine the values of $D_{\ti+1,\ka}$ for $\ti\geq 0$. 
For the random perturbations $\epsilon_{\ti+1,\ka}$, we generate a single fixed sample of values $\{\hat{\epsilon}^\ell_{\ka}: \ell=1,\ldots,m, \ka \in \KA \}$ according to the standard normal distribution. Then, in each iteration of simulating demands following the AR process, we choose a sample from this fixed set uniformly at random.

The algorithms are implemented in Python and MILP/LP models are solved using IBM ILOG CPLEX version 12.8. All the experiments are performed on a 2.4 GHz Intel Gold processor with only one thread on the Béluga, Digital Research Alliance of Canada computing grid.

\subsection{Methodology evaluation}
\label{ss:method}
In this section, we test the efficiency of the two methods for solving the two-stage chance-constrained model, the extensive form described in Section \ref{sec:ext} and the proposed branch-and-cut (B\&C) algorithm described in Section \ref{sec:bnc}. To this end, we generate one instance  of each of the following parameter combinations: $\Ti \in \{  6, 8 , 10\}$, $\Ka = 10$,
$\eta =  0.2 $,
$\tau = 0.5 $,
$\rho = 0.1$,
$ TBO \in \{1, 2\} $,
$ \alpha \in \{80\%, 90\%, 95\%, 99\% \}$, and $|\Omega| \in \{100,200,300,500,1000\}$. Thus, we have 120 instances in total.
The initial state is set by running the simulation described in the next section through its warm-up period using a fixed policy, and then using the initial state for the five iterations using the two policies. We emphasize that the policy used in the simulation is used just to generate the initial state for the next five iterations. Given one such fixed instance, we then solve it with the two different methods to compare the computational performance of the methods. We set a time limit of 7200 seconds to solve each of these instances.

  
We analyze the performance of the two methods using three measures: the average CPU time in seconds (Time), the average optimality gap after the time limit is reached (Gap), and the percentage of instances that were solved to optimality within the time limit (\% OPT).

\begin{table}[ht]
\small
\centering
\caption{Comparison of methodologies to solve the two-stage chance-constrained model.}
\label{tab:MethodologyCompare}
\renewcommand{\arraystretch}{0.8}
\begin{tabular}{lrccrcc}
\toprule
  & \multicolumn{3}{c}{ B\&C} &  \multicolumn{3}{c}{Extensive form}
  \\
\cmidrule(lr){2-4}
\cmidrule(lr){5-7}
& Time       & Gap (\%)       & \% Opt      & Time    & Gap(\%)     & \% Opt   \\
\midrule
$\Em$  &&&&&&\\
100	&	10.3	&	0.0	&	100	&	74.6	&	0.0	&	100	\\
200	&	34.6	&	0.0	&	100	&	400.8	&	0.0	&	100	\\
300	&	60.6	&	0.0	&	100	&		1152.5	&	0.0	&	97	\\
500	&	206.1	&	0.0	&	100	&	3356.3	&	0.4	&	80	\\
1000	&	990.1	&	0.0	&	98	&	6170.3	&	1.7	&	26	\\
\midrule
$\alpha$ (\%) &&&&&&\\
80&	417.2	&	0.0	&	99	&		2486.6	&	0.4	&	78	\\
90&	299.5	&	0.0	&	100	&	2690.7	&	0.6	&	76	\\
95&	202.2	&	0.0	&	99	&	2166.3	&	0.5	&	79	\\
99&	122.5	&	0.0	&	100	&	1579.9	&	0.4	&	89	\\
\midrule
$\Ti$ &&&&&&\\
6	&	131.4	&	0.0	&	100	&		1780.3	&	0.4	&	84	\\
8	&	261.8	&	0.0	&	99	&	2040.1	&	0.3	&	84	\\
10	&	387.9	&	0.0	&	99	&		2872.3	&	0.8	&	73	\\
\midrule
Average &	260.4	&	0.0	&	100	&		2230.9	&	0.5	&	80	\\
\bottomrule
\end{tabular}

\end{table}

The results are given in Table \ref{tab:MethodologyCompare}, where each row presents results averaged over all instances sharing a particular parameter level as given in the first column. For example, the first row of data presents aggregate results over all instances with $|\Omega|=100$, and the first row of data in $\alpha$ section presents aggregate results over all instances with $\alpha=80\%$. From this table we see that the branch-and-cut method successfully solves nearly all instances within the time limit, and in significantly less time than using the extensive form, and that this result is consistent across all ranges of parameters. Most significantly, we observe that with the branch-and-cut method we are able to solve instances of varying size in terms of number of time periods and number of scenarios used to approximate the chance constraint.
The results also indicate that instances with higher service level can be solved faster by both methods. Finally, as expected we observe that the solution time increases with the number of time periods and with the number of scenarios used to approximate the distribution of product demands.


\subsection{Policy evaluation}

\subsubsection{Evaluation via simulation}
\label{Sec:Rolling}

Recall that the setting for the problem we study is an infinite-horizon problem in which decisions are repeatedly made over time, and the proposed decision policies are based on solving a finite-horizon problem to be used in a rolling-horizon framework. 
That is, in each period a model with $\Ti$ planning periods is solved, and only the decisions corresponding to the first period are implemented. Based on these decisions and the observed demand, the state of the system is updated and the next $\Ti$-period model is solved, and the process repeats. 

\begin{algorithm}[ht]
\SetAlgoLined
\linespread{1.15}\selectfont
{OUTPUT: The confidence intervals on the total cost and the service level}\\
INPUT: \
 Demand simulation over $\Ti_{Sim}$ periods, Number of warm-up periods $T_{Warm}$, Production policy\\

 $\tAct =1,\hat{\Vi}_{0} = \boldsymbol{0},\hat{\Bi}_{0} = \boldsymbol{0}, \mathcal{O} = \emptyset, \mathcal{Z} = \emptyset, \hat{D}_{1} = \boldsymbol{0}$ \\
 \While{$\tAct \leq \Ti_{Sim}$ }{
 Solve the stock-out determination model~(\ref{Currentstage}) using  $(\hat{D}_{\tAct},\hat{\Bi}_{\tAct-1},\hat{v}_{\tAct-1})$ as the input
 $(\hat{D}_1,\hat{\Bi}_0,\hat{v}_0)$ and let $\Bi^{*}$ be its optimal backlog solution.
 \\
 Let $\hat{\KA} = \{ k \in \KA : \Bi^{*}_{\ka} = 0 \}$

  Solve either \eqref{mod:Det} or a sample-approximation of \eqref{mod:ChanceConstraintpolicy}, based on selected policy, using 
  $(\hat{D}_{\tAct},\hat{\Bi}_{\tAct-1},\hat{v}_{\tAct-1})$ as the input $(\hat{D}_1,\hat{\Bi}_0,\hat{v}_0)$ and the computed set $\hat{\KA}$.\\
 Let $\hat{x}_{\tAct}, \hat{y}_{\tAct}, \hat{\Es}_{\tAct}, \hat{\Bi}_{\tAct}, \hat{i}_{\tAct},\hat{v}_{\tAct}$ be the first-period components of the optimal solution, and let the $Obj_{\tAct}$ be the total cost in period $\tAct$ based on this solution. \\ 
  \eIf{$\exists \ka \in \KA$ with $\Bi^{*}_{\tAct \ka} > 0$} {\vskip 0.1cm    
  $Z_{\tAct} = 1$ }
  {
            $Z_{\tAct} = 0$
        }
  \If{$\tAct \geq T_{Warm}$}{
  Add $Obj_{\tAct}$ and $Z_{\tAct}$ to the set $\mathcal{O}$ and $\mathcal{Z}$, respectively \\
  }
  \For{$\ka \in \KA$}{
   \If{$\tAct =1$}{ $\hat{D}_{\tAct k} \gets C/(1-AR_1)$ 
   }
   Obtain a realization $\hat{\epsilon}$ of $\epsilon_{\tAct+1,\ka}$ \\
    $\hat{D}_{\tAct+1,\ka} \gets  C + AR_1 \hat{D}_{\tAct k} + AR_2 \hat{\epsilon}$
  }
  $\tAct \gets \tAct+1$ }
  Build confidence intervals for the cost and service level using $\mathcal{O}$ and $\mathcal{Z}$, respectively.
  \caption{Rolling-horizon implementation}
  \label{Al:PolicyEvaluation}
\end{algorithm}

We implement a steady-state simulation to test different policies, as described in Algorithm \ref{Al:PolicyEvaluation}. At each time period $\tAct$ we first execute the stock-out determination model and then solve a finite-horizon model depending on the selected policy. 
Specifically, we test three different policies:
\begin{itemize}
\item {\it Average}: Based on solving the deterministic approximation \eqref{mod:Det}, where we use the expected value of future demands as the demands in periods $2,\ldots,T$ in \eqref{mod:Det}, where the expected values are conditional on the current observed demands.
\item {\it Quantile}: Based on solving the deterministic approximation \eqref{mod:Det}, but with the $\alpha$ quantile of the random demand used as the demand for each product in period 2, and the expected values of future demands are used as the demands in periods $3,\ldots,T$.
\item {\it CC}: Based on solving a sample-average approximation of the two-stage chance-constrained model \eqref{mod:ChanceConstraintpolicy}.
\end{itemize}
For all policies we use a look-ahead horizon of $T=6$ periods, which was determined based on preliminary experiments that indicated using more periods did not appear to yield better results. For the CC policy, we use a sample size of 100 scenarios for the sample average approximation. We use this relatively small number of scenarios to ease the computational burden of the experiments, since we must solve this model in each of the (over 4000) periods of the simulation run for each test instance. We emphasize that when using this policy in practice it would only be necessary to solve a single model in each period, and the results from Section \ref{ss:method} demonstrate that it would then be feasible to use significantly more scenarios (and a longer time horizon to look ahead) in case that is necessary to yield a better approximation.

We run the simulation for $\Ti_{Sim}=4010$ time periods, and ignore the first $\Ti_{Warm}=10$ time periods as a warm-up phase when computing estimates of the average cost and service level. In each time period after the warm-up phase, we record the actual cost (sum of setup costs, holding costs, and substitution costs) in that time period and an indicator of whether or not there was backlog in any product during that time period.  Recall that we do not include production costs, as the long-run average of the total number of products made per period is the same for all policies, and differences in production costs that are incurred due to substitution are included in the substitution cost.
For calculating the confidence intervals on these measures, we use batch-means estimation, with 160 batches of 25 time periods each.






\subsubsection{Policy comparison}

We first compare the three policies against each other.
This comparison is based on the estimated average cost per period and the estimated service level, using the procedure explained in Section~\ref{Sec:Rolling}. Table \ref{PolicyComp} compares the three policies using these measures and their 95\% confidence interval at two different levels of the TBO parameter and four different service level targets. 
Among the three policies the CC policy is the only policy which respects the service level target in all the instances. In all the instances with acceptable service level the CC policy has the lowest cost. Among the three policies, the average policy consistently provides service level below the target.  We thus do not report the performance of the average policy in the following experiments. The quantile policy, on the other hand, demonstrates better potential for meeting the service level targets.

\newcolumntype{L}{>{$}l<{$}}
\newcolumntype{C}{>{$}c<{$}}
\newcolumntype{R}{>{$}r<{$}}
\newcommand{\nm}[1]{\textnormal{#1}}

\begin{table} [htpb]
\small
\centering
\caption{Policy comparison based on total cost and service level}\label{PolicyComp}
\small
\renewcommand{\arraystretch}{1}
\begin{tabular}{RRRRRRRR}
\toprule
\multicolumn{1}{R}{$TBO$} &
\multicolumn{1}{R}{\alpha \ (\%)} &
\multicolumn{3}{c}{Total cost}    &
\multicolumn{3}{c}{Service level (\%)}    \\ 
\cmidrule(lr){3-5}
\cmidrule(lr){6-8}

&&
\multicolumn{1}{r}{Average } &
\multicolumn{1}{r}{Quantile }     &
\multicolumn{1}{r}{CC } &
\multicolumn{1}{r}{Average } &
\multicolumn{1}{r}{Quantile }     &
\multicolumn{1}{r}{CC}  \\
\midrule

\multicolumn{1}{c}{1}	&\multicolumn{1}{c}{	80}	&	74.4	\pm	0.4	&	66.4	\pm	0.2	&	66.7	\pm	0.2	&	21.4	\pm	1.4	&	76.5	\pm	1.4	&	84.9	\pm	1.3	\\
	&	\multicolumn{1}{c}{90}	&	74.4	\pm	0.4	&	68.2	\pm	0.1	&	67.1	\pm	0.2	&	21.4	\pm	1.4	&	90.6	\pm	1.1	&	90.3	\pm	1.1	\\
	&\multicolumn{1}{c}{	95}	&	74.4	\pm	0.4	&	71.0	\pm	0.1	&	67.6	\pm	0.2	&	21.4	\pm	1.4	&	94.1	\pm	0.9	&	95.3	\pm	0.7	\\ \vspace{0.2cm}
	&\multicolumn{1}{c}{	99}	&	74.4	\pm	0.4	&	76.1	\pm	0.1	&	69.1	\pm	0.2	&	21.4	\pm	1.4	&	99.1	\pm	0.3	&	99.1	\pm	0.3	\\ 
\multicolumn{1}{c}{2}	&\multicolumn{1}{c}{	80}	&	204.3	\pm	0.6	&	204.3	\pm	0.5	&	191.2	\pm	0.6	&	78.1	\pm	2.7	&	93.2	\pm	1.2	&	98.7	\pm	0.4	\\
	& \multicolumn{1}{c}{	90}	&	204.3	\pm	0.6	&	207.2	\pm	0.4	&	192.4	\pm	0.6	&	78.1	\pm	2.7	&	97.4	\pm	0.6	&	99.5	\pm	0.3	\\
	&\multicolumn{1}{c}{	95}	&	204.3	\pm	0.6	&	210.0	\pm	0.4	&	193.5	\pm	0.6	&	78.1	\pm	2.7	&	98.5	\pm	0.5	&	99.8	\pm	0.2	\\
	&\multicolumn{1}{c}{	99	}&	204.3	\pm	0.6	&	215.3	\pm	0.4	&	195.4	\pm	0.6	&	78.1 \pm	2.7	&	99.7	\pm	0.2	&	100.0	\pm	0.0	\\

\bottomrule
\end{tabular}
\end{table}

\begin{figure}[b]
    \centering
    \subfloat[\centering Service level]{{\includegraphics[width=9cm]{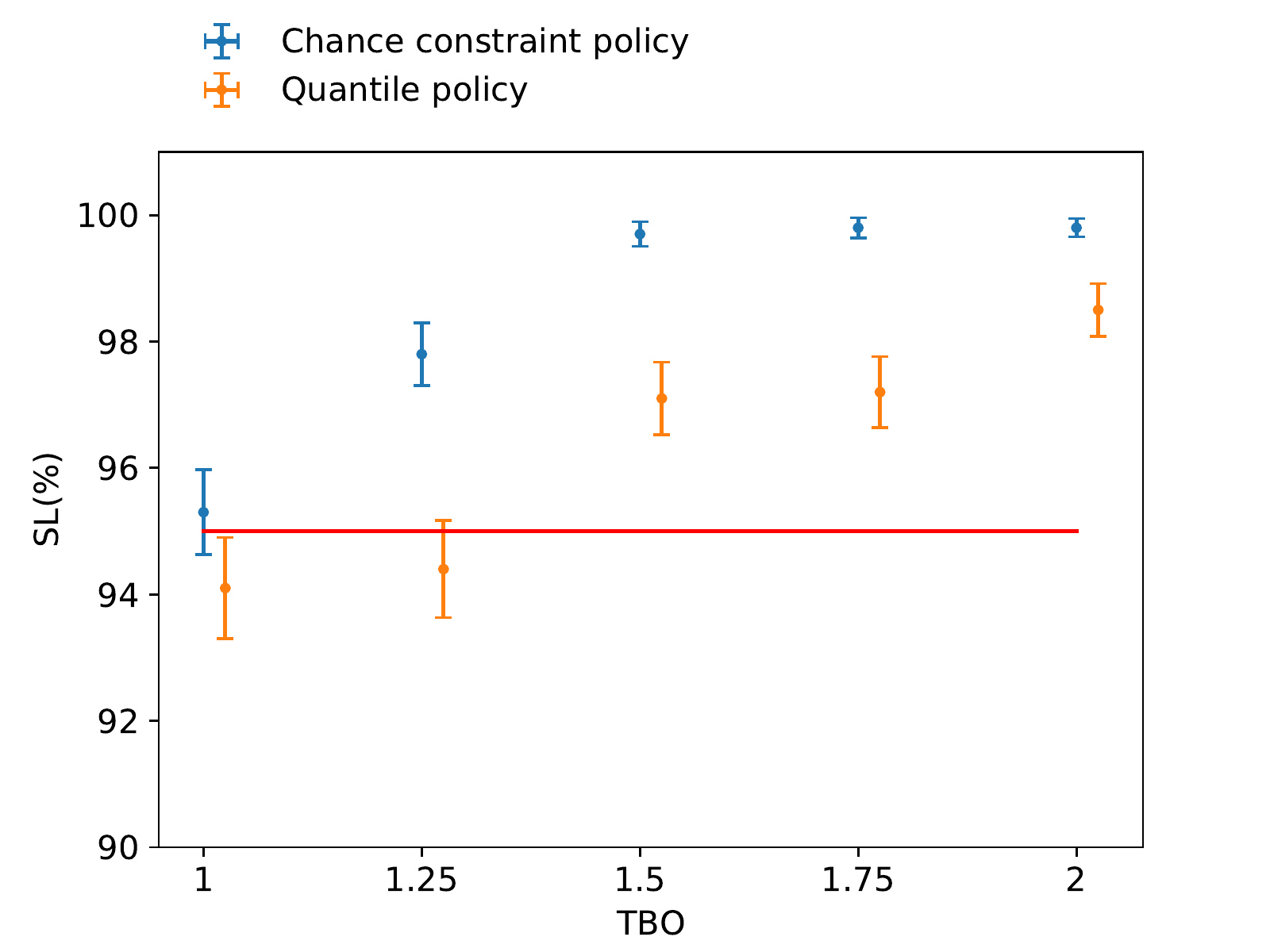} }}%
    \subfloat[\centering Total cost ]{{\includegraphics[width=9cm]{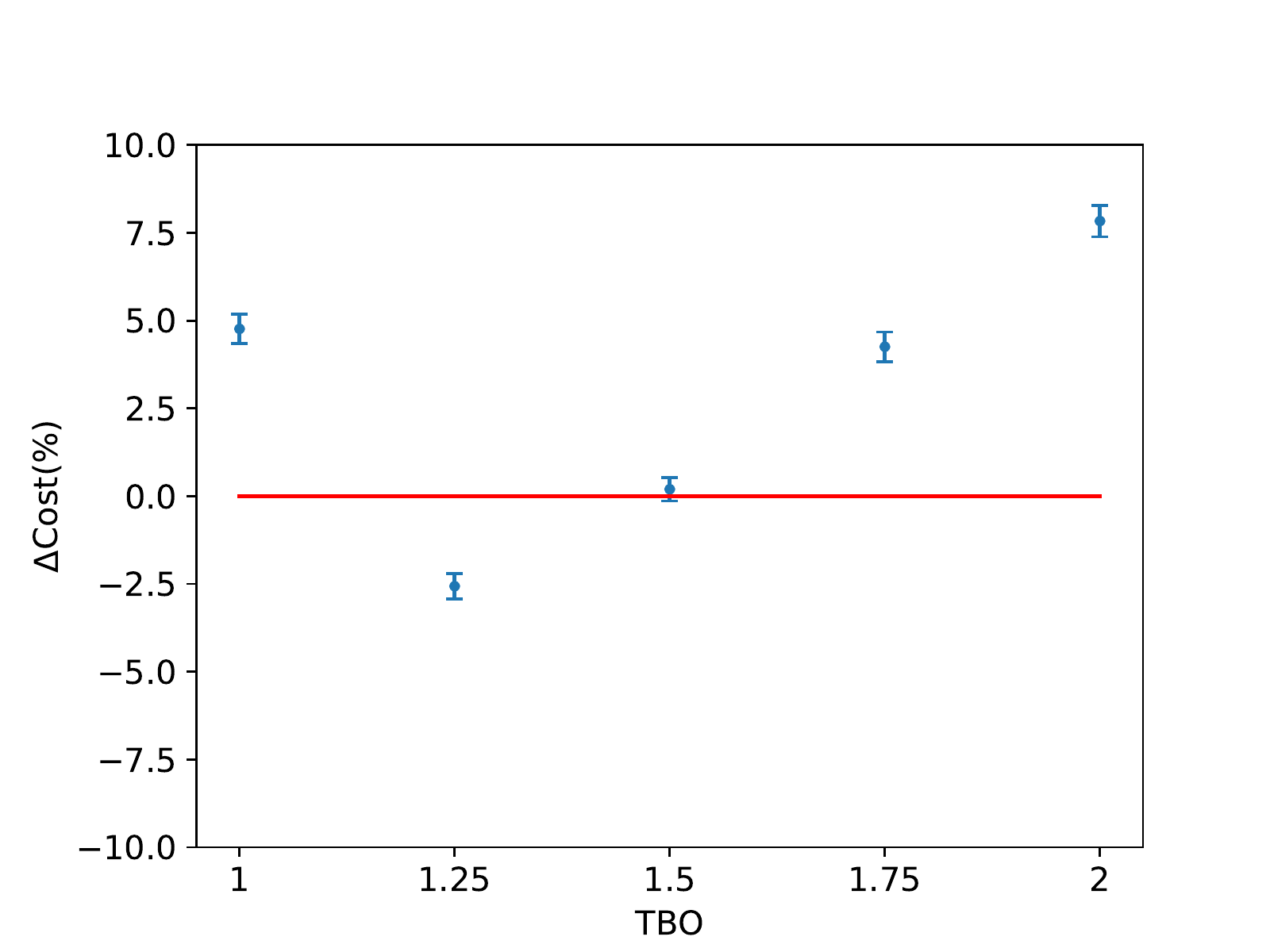} }}%
    \caption{Comparison based on TBO}%
    \label{fig:TBOComp}%
\end{figure}

The rest of this section is dedicated to the comparison between the CC and quantile policies using the additional instances presented in Table~\ref{tab:BaseSensitivity}. To this end, we use two measures, the joint service level and the relative cost change, $\Delta$Cost, which is defined as: 
\begin{alignat}{2}
  & \Delta \text{Cost} (\%) = \dfrac{\text{Total Cost}_{\text{Quantile}} - \text{Total Cost} _{\text{CC}}}{\text{Total Cost}_\text{{Quantile}}} \times 100& \label{eq:ِDeltaCost} 
 \end{alignat}
Positive $\Delta$Cost means that the CC policy had lower costs than the Quantile policy.
Figure~\ref{fig:TBOComp} shows the comparison of the quantile policy and the CC policy under different values of TBO. Figure~\ref{fig:TBOComp}-(a) shows the service level, labeled as SL, for each of the policies at different values of TBO and Figure~\ref{fig:TBOComp}-(b) illustrate the $\Delta$Cost for each value of TBO. In these and all following figures in this section, the point estimates of the quantities are displayed with a point and the whiskers represent the 95\% confidence interval of the estimated quantity.
In all cases, the CC policy has a better performance in terms of service level. The CC policy has a lower cost in all cases in which both policies have an acceptable service level. When TBO is more than 1 the obtained service level exceeds the target.  This is caused by a combination of the impact of higher setup cost when TBO is larger than 1 and the use of the stock-out determination step to enforce that backlog is positive only when necessary. Specifically, when the setup costs are higher, it is generally optimal to place orders less frequently, and thus more inventory is carried on average. When there is more inventory on-hand, the stock-out determination step will usually find that it is possible to avoid any backlog. 

\begin{figure}[htbp]%
    \centering
    \subfloat[\centering Service level]{{\includegraphics[width=8.6cm]{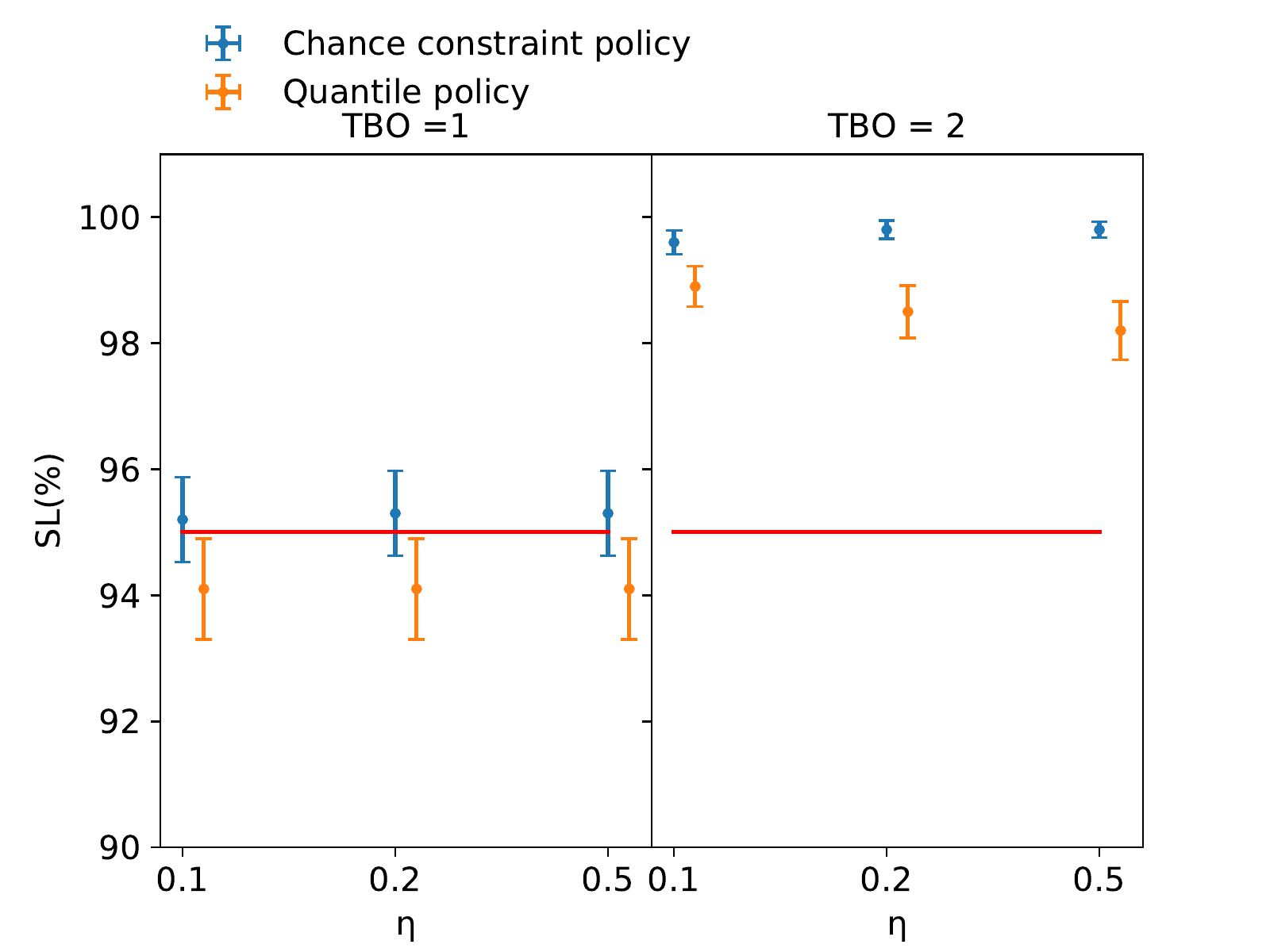} }}%
    \subfloat[\centering Total cost ]{{\includegraphics[width=8.6cm]{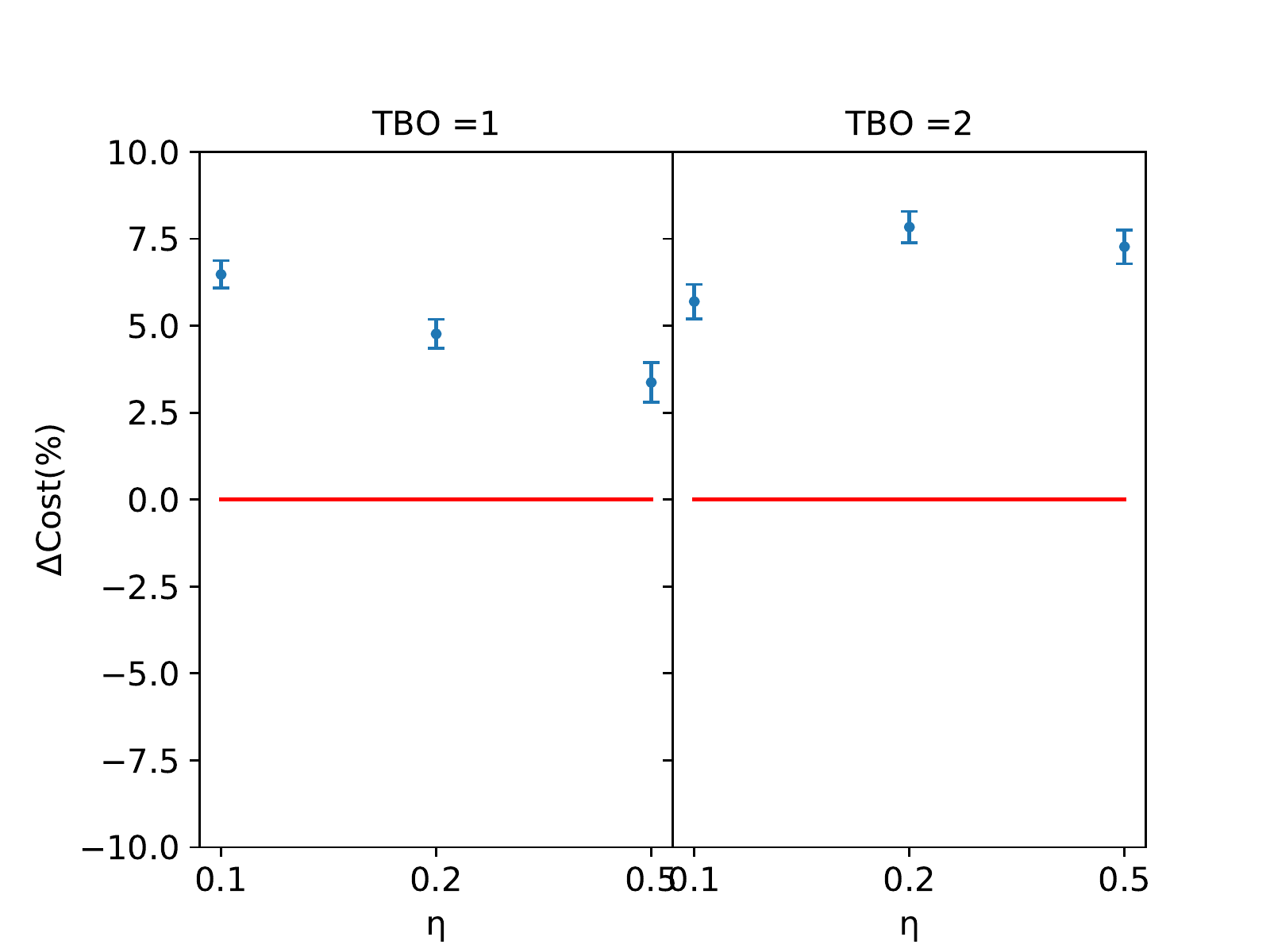} }}%
    \caption{Comparison based on $\eta$}%
    \label{fig:ETAComp}%
\end{figure}

Figure~\ref{fig:ETAComp} shows the comparison based on different values of the production cost variability parameter $\eta$ under two different values of TBO, 1 and 2.  Higher $\eta$ means higher variability in the production costs.  When TBO is equal to 1, the quantile policy service level is slightly lower than the target service level. In all cases, the CC policy has a better performance in terms of service level and total cost. 

\begin{figure}[htbp]
    \centering
    \subfloat[\centering Service level]{{\includegraphics[width=8.6cm]{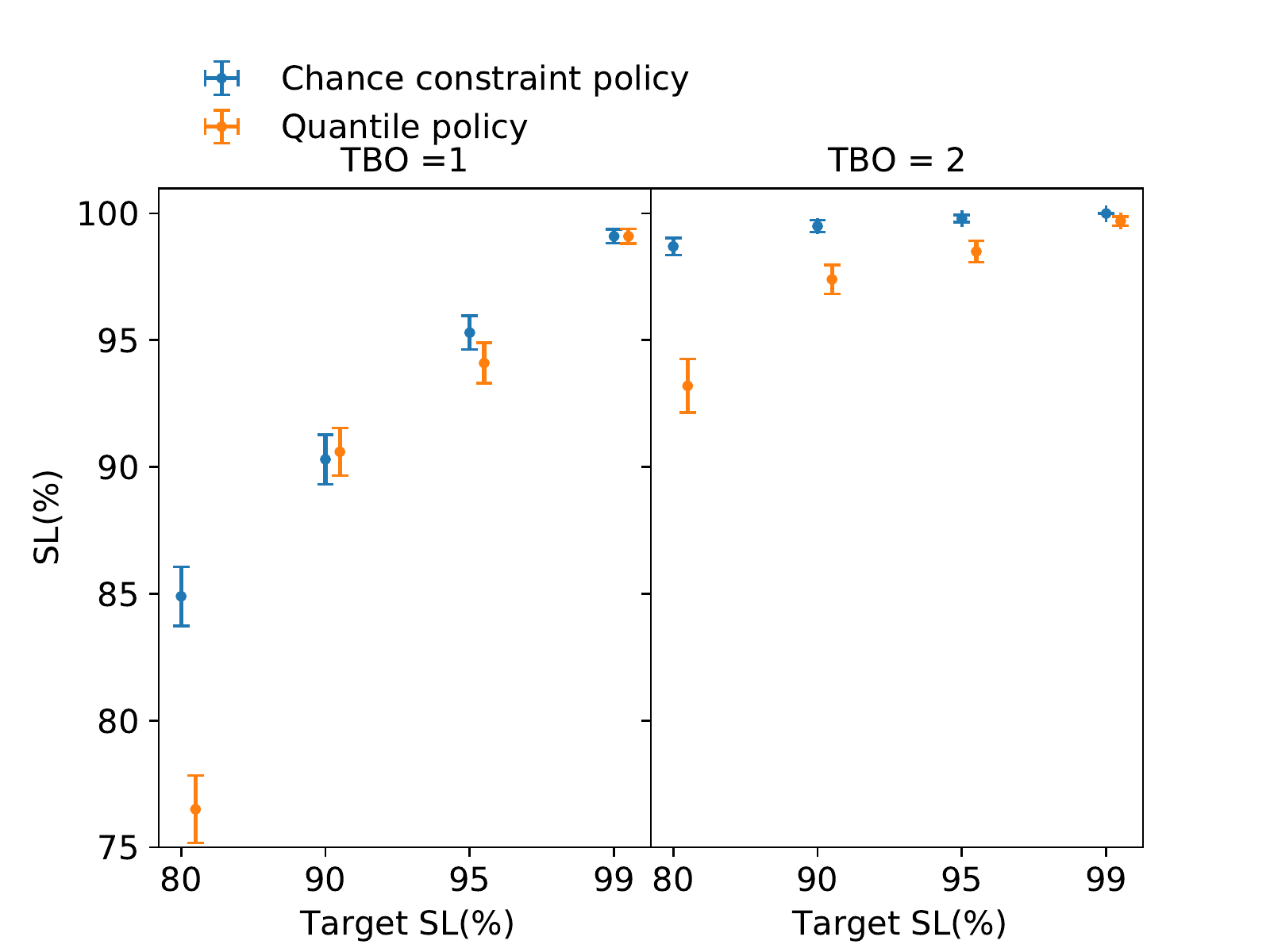} }}%
    \subfloat[\centering Total cost ]{{\includegraphics[width=8.6cm]{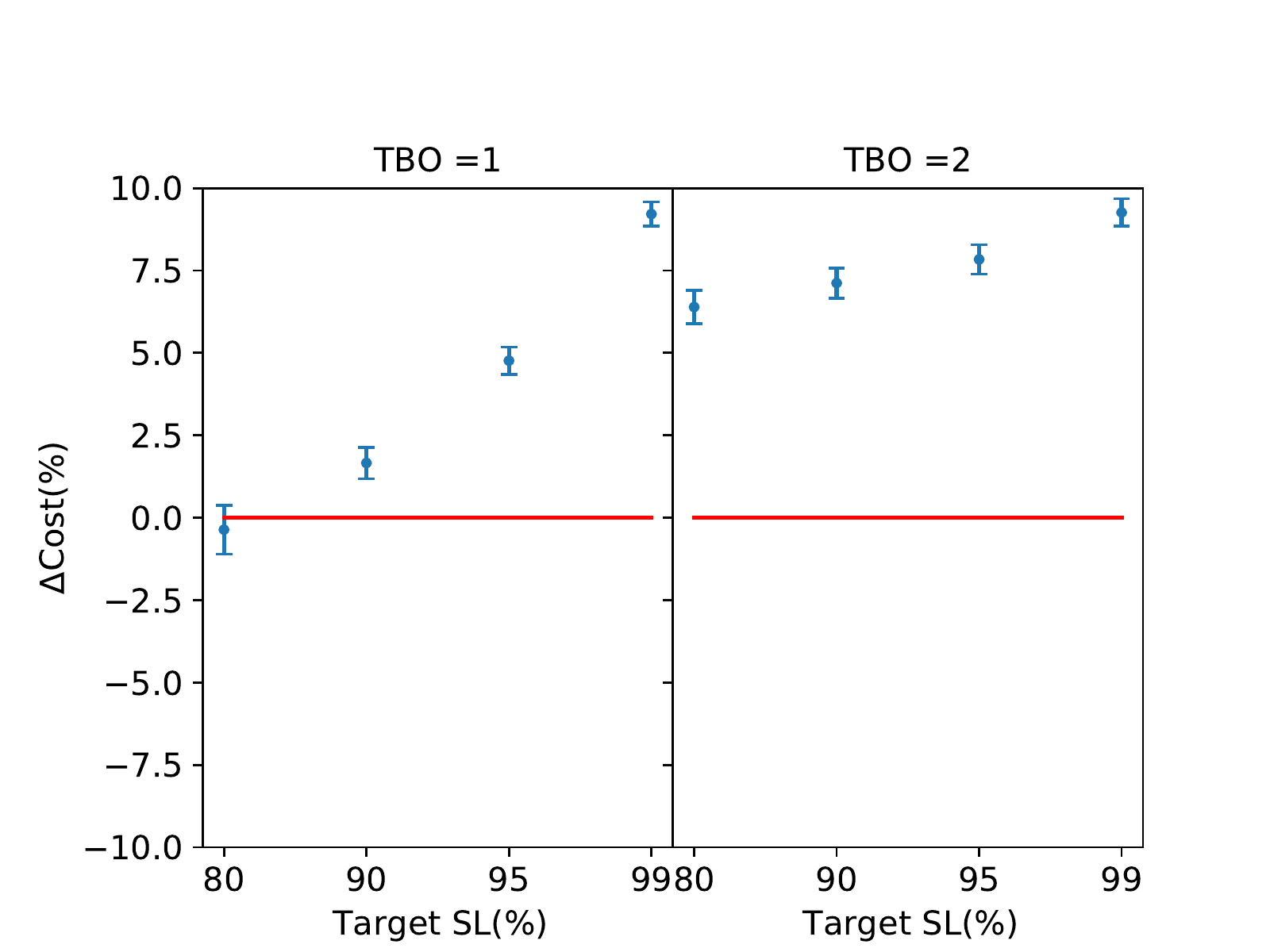} }}%
    \caption{Comparison based on the target service level}%
    \label{fig:SLComp}%
\end{figure}

Figure~\ref{fig:SLComp} shows the comparison based on different service level targets under two different values of TBO, 1 and 2. In all cases, the CC policy respects the service level target and in cases where both policies have an acceptable service level, the CC policy has better performance in terms of the total cost. We observe that when the service level increases, the relative performance of the CC policy against the quantile policy improves. 

\begin{figure}[ht]
\begin{center}
\includegraphics[scale=0.6]{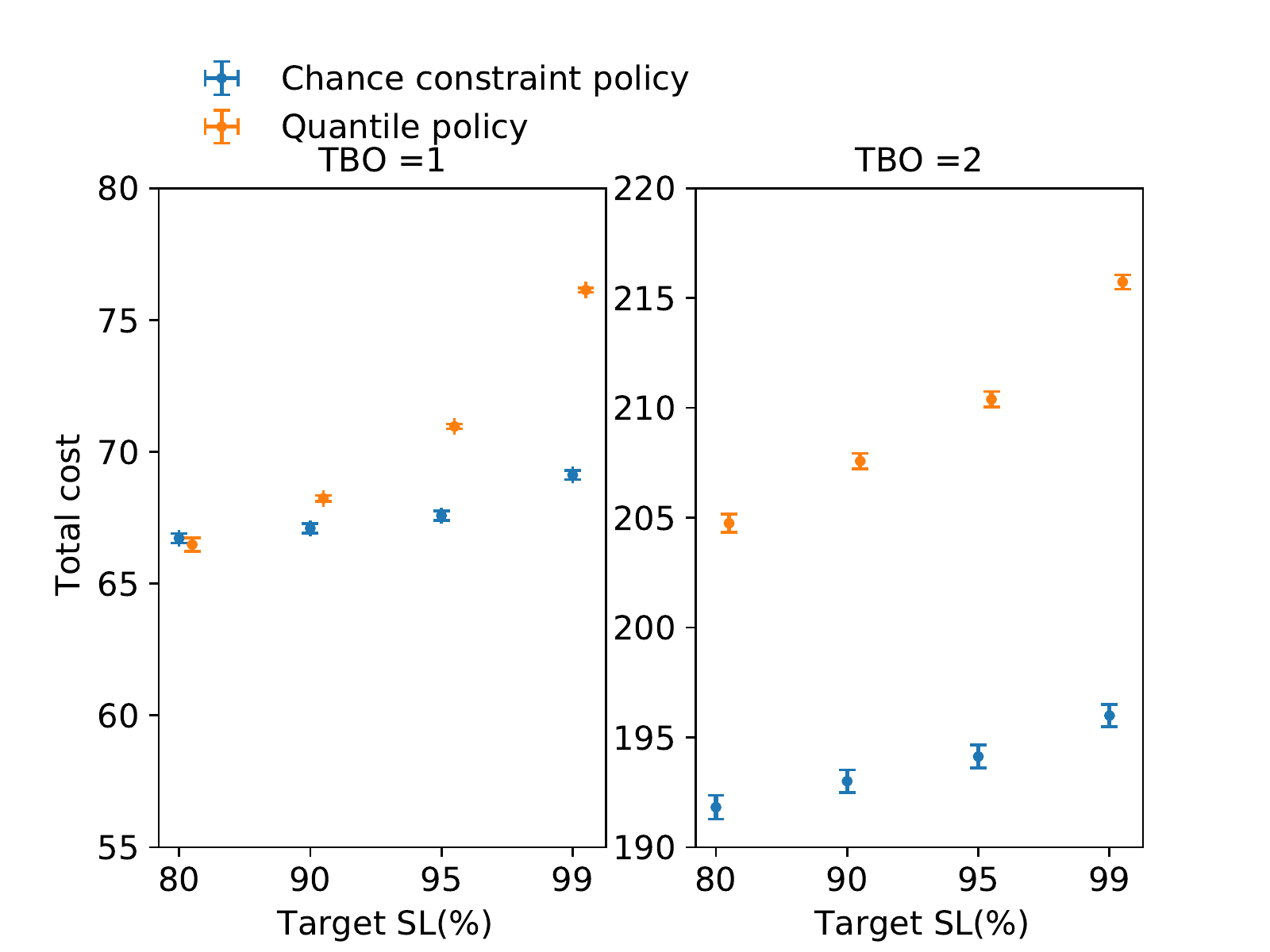}
\caption{Total cost trend comparison based on $\alpha$} 
\label{fig:SLTotalCost}
\end{center}
\end{figure}

Figure~\ref{fig:SLTotalCost} is complementary to Figure~\ref{fig:SLComp} and illustrates the trend of the total cost for different values of service level. As can be seen in this figure, the total cost of the quantile policy increases significantly with an increase in the target service level, whereas with the CC policy a higher service level can be achieved with significantly less increase in cost.

\begin{figure}[ht]
    \centering
    \subfloat[\centering Service level]{{\includegraphics[width=8.2cm]{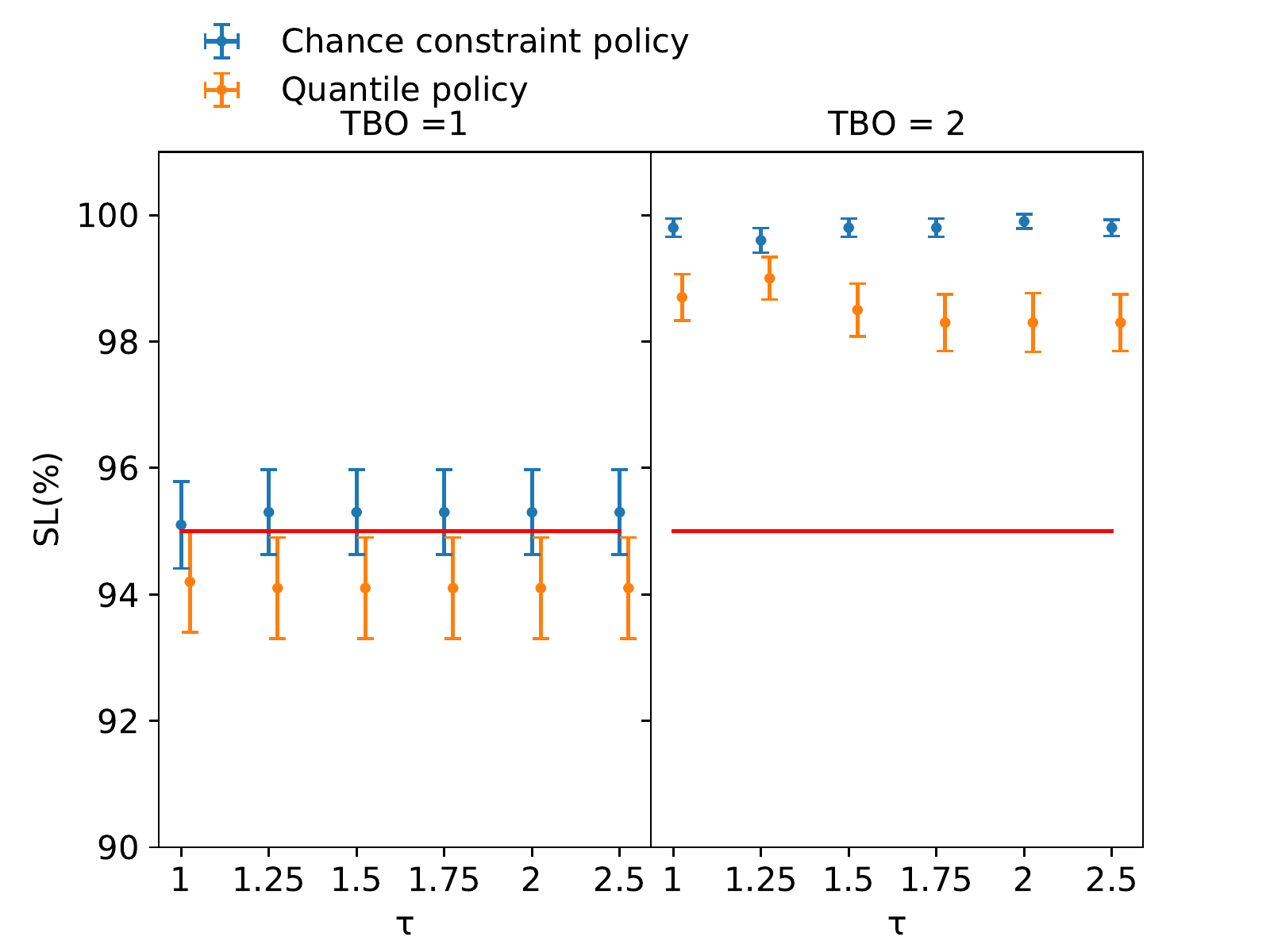} }}%
    \subfloat[\centering Total cost ]{{\includegraphics[width=9cm]{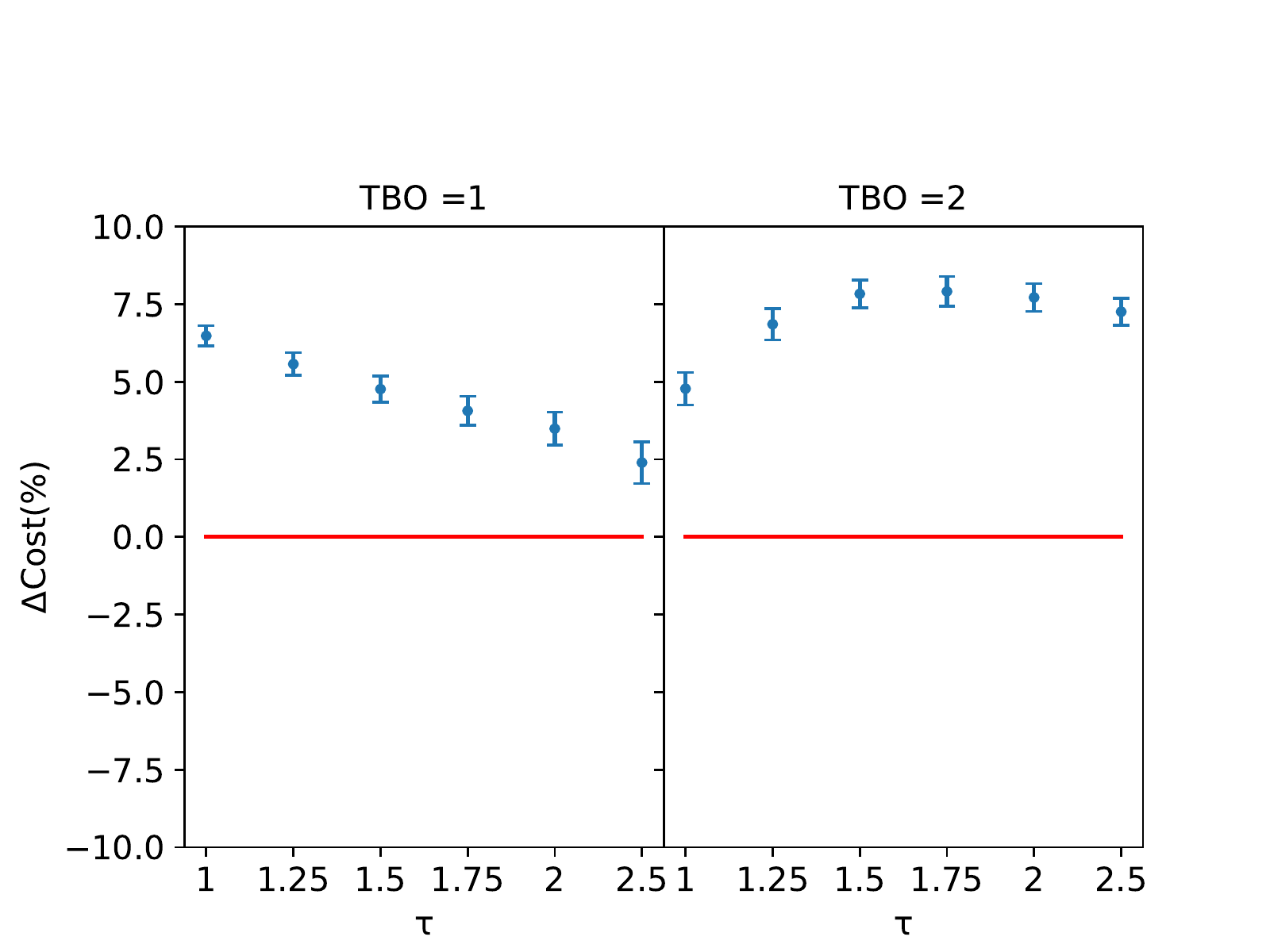} }}%
    \caption{Comparison based on $\tau$}%
    \label{fig:TIComp}%
\end{figure}

Figure~\ref{fig:TIComp} shows a similar comparison based for varying values of the parameter $\tau$,  which impacts the substitution cost ($\tau=1$ is the minimum case where substitution cost is just equal to the difference in production costs, whereas $\tau >1$ adds a higher penalty for substitution). We see that the CC policy yields higher service levels and lower costs than the quantile policy over all tested values of $\tau$.

We can conclude that although the quantile policy has a reasonable performance in general in terms of meeting the service level target,  the proposed CC policy consistently achieves both higher service levels and lower costs than the quantile policy.

\subsection{The importance of the stock-out determination step}

    When TBO is greater than 1 the service level obtained with the CC policy exceeds the target  (see Figure~\ref{fig:TBOComp}-(a)). This is because setup costs are higher when TBO is greater than 1, and hence when production occurs, the production quantities are higher to save on setup costs. This leads to higher inventory levels on average, and hence it is frequently possible to avoid having any backlog in a period. However, this raises the question of whether average costs could be reduced further if we did not use the stock-out determination model to avoiding backlog whenever possible. We thus conduct an experiment to estimate the service level and total cost when the stock-out determination step is skipped. Specifically, in this case we never enforce that the backlog variables equal to zero when solving the chance-constrained model \eqref{mod:ChanceConstraintpolicy}.
    
    \begin{figure} [ht]
    \centering
    \subfloat[\centering Service level]{{\includegraphics[width=8.6cm]{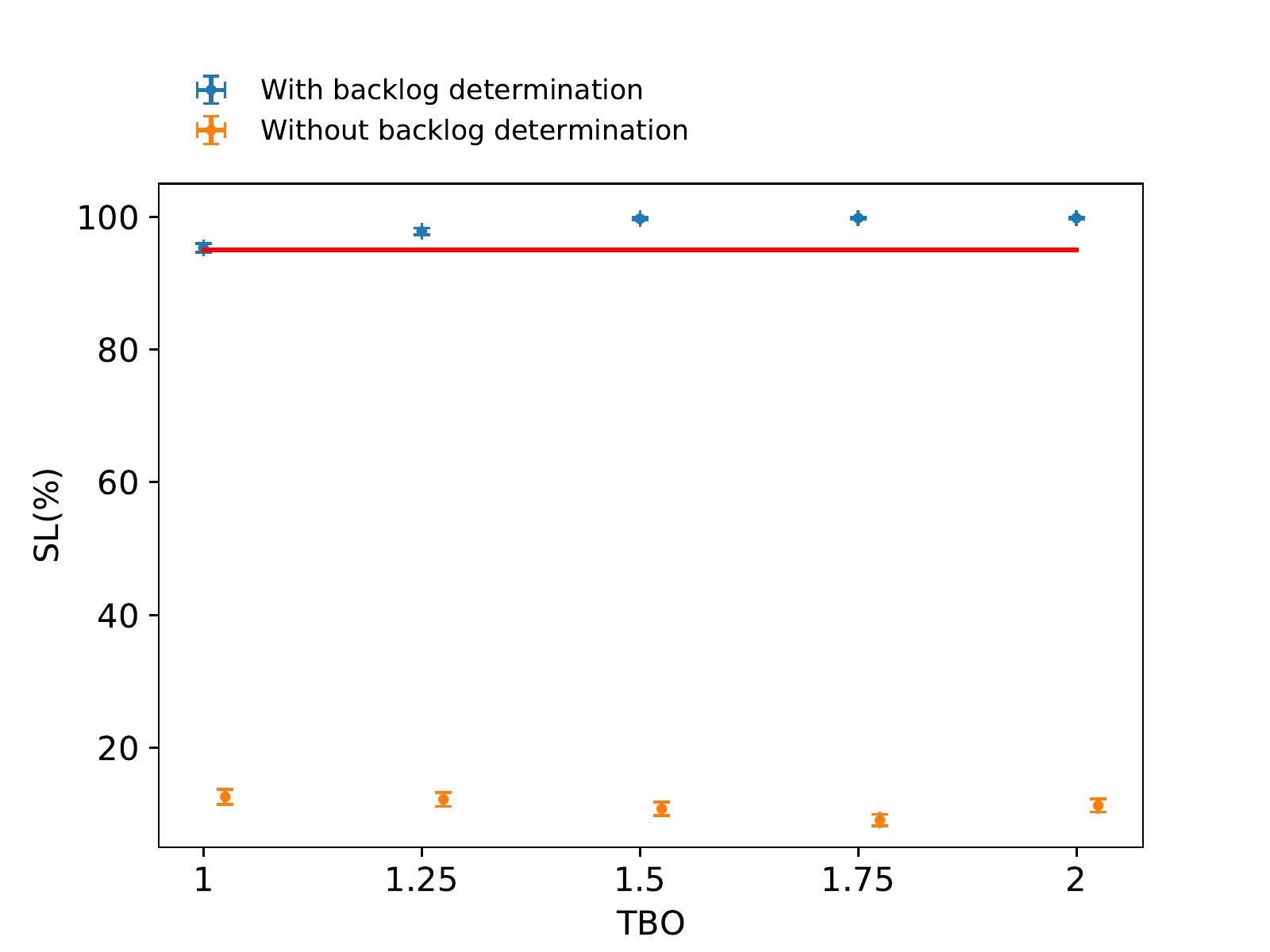} }}%
    \subfloat[\centering Total cost ]{{\includegraphics[width=8.6cm]{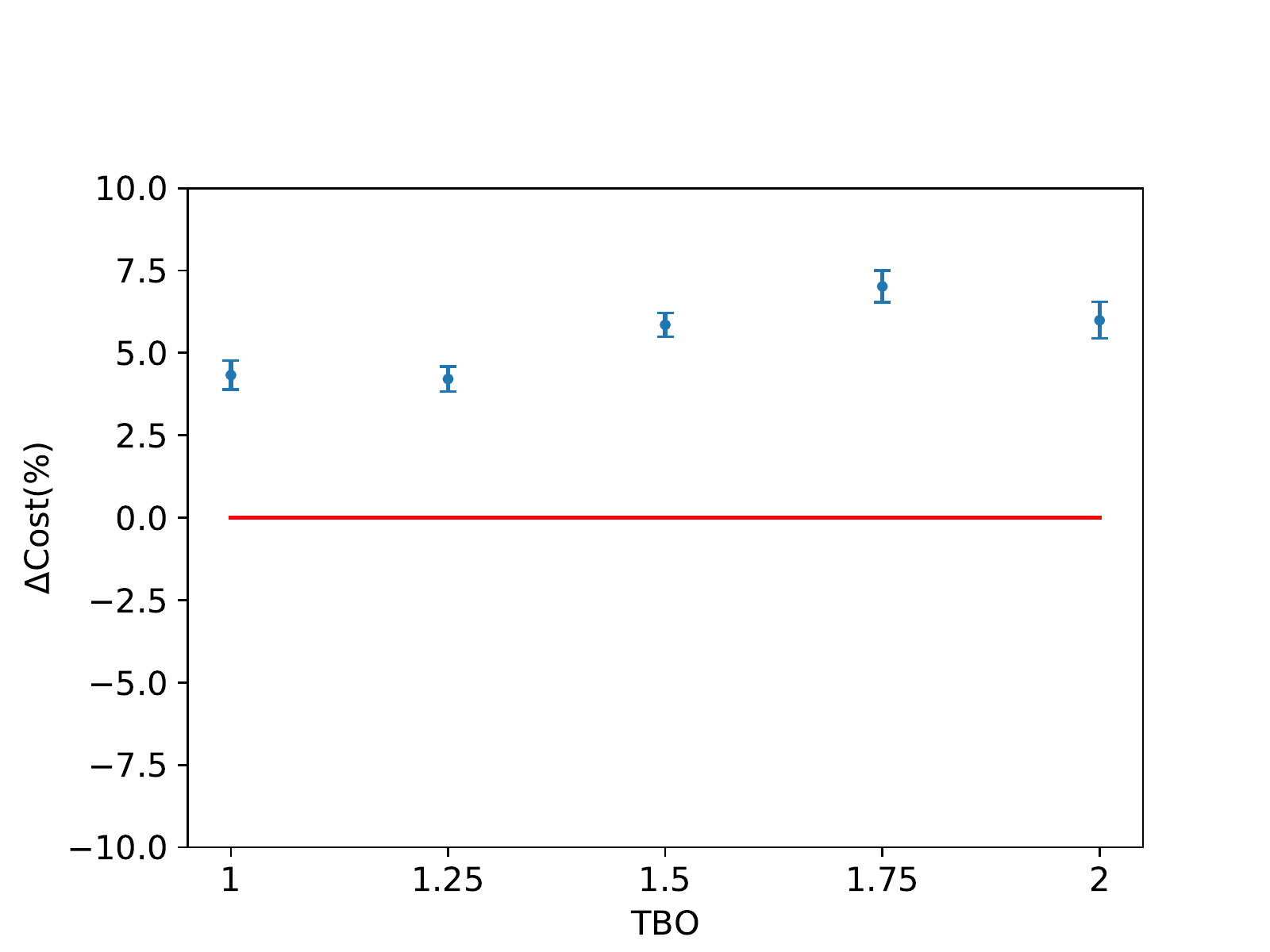} }}%
    \caption{The necessity of the stock-out determination step.}%
    \label{fig:BacklogDetermination}%
\end{figure}
    Figure~\ref{fig:BacklogDetermination}-(a) presents the service level obtained with and without the stock-out determination step and
    Figure~\ref{fig:BacklogDetermination}-(b) presents the cost decrease that is obtained when the stock-out determination step is skipped. We observe that skipping the stock-out determination step does lead to cost reductions, but the magnitude of the reductions is modest. On the other hand, without the stock-out determination step the service level falls below 20\%. This gap indicates that there are many periods in which it is possible to meet all customer demands, but solving the chance-constrained model \eqref{mod:ChanceConstraintpolicy} without any constraints requiring this consistently leads to solutions in which demands are backlogged. This illustrates the need to have some mechanism that assures demands are fulfilled in the current period when possible. While the stock-out determination step is not the only possibility for achieving this, this experiment suggest that it is reasonably effective, as the cost increase is modest even when compared to the extreme alternative of ignoring current period demands altogether. 

The stock-out determination step can be interpreted as allowing backlogging only when absolutely necessary. While we do not pursue this here, a conceptually simple modification to this policy would be to allow backlogging when the cost of meeting all current demands (e.g., via substitutions) exceeds some fixed threshold. This threshold would need to be tuned so that the service level target is satisfied. This may allow a reduction in average cost by reducing the degree to which the achieved service level exceeds the target.

\subsubsection{Effect of substitution}

We next investigate the extent to which substitution allows achieving service level targets at reduced costs. We also explore the relative benefits from allowing a wider range of products to be substituted for each other.  To this end, we evaluate the service level and average cost using the CC policy and three levels of substitution: (1) no substitution allowed, (2) partial substitution, which corresponds to our base case in which product
$k$ can be used to meet demand of product $j$ if $j \geq k$ (so $k$ is a higher quality product) and $k \geq j-3$, and (3) full substitution, in which product $k$ can be used to meet demand of any product $j \geq k$.

\begin{figure}[ht]
\begin{center}
\includegraphics[scale=0.6]{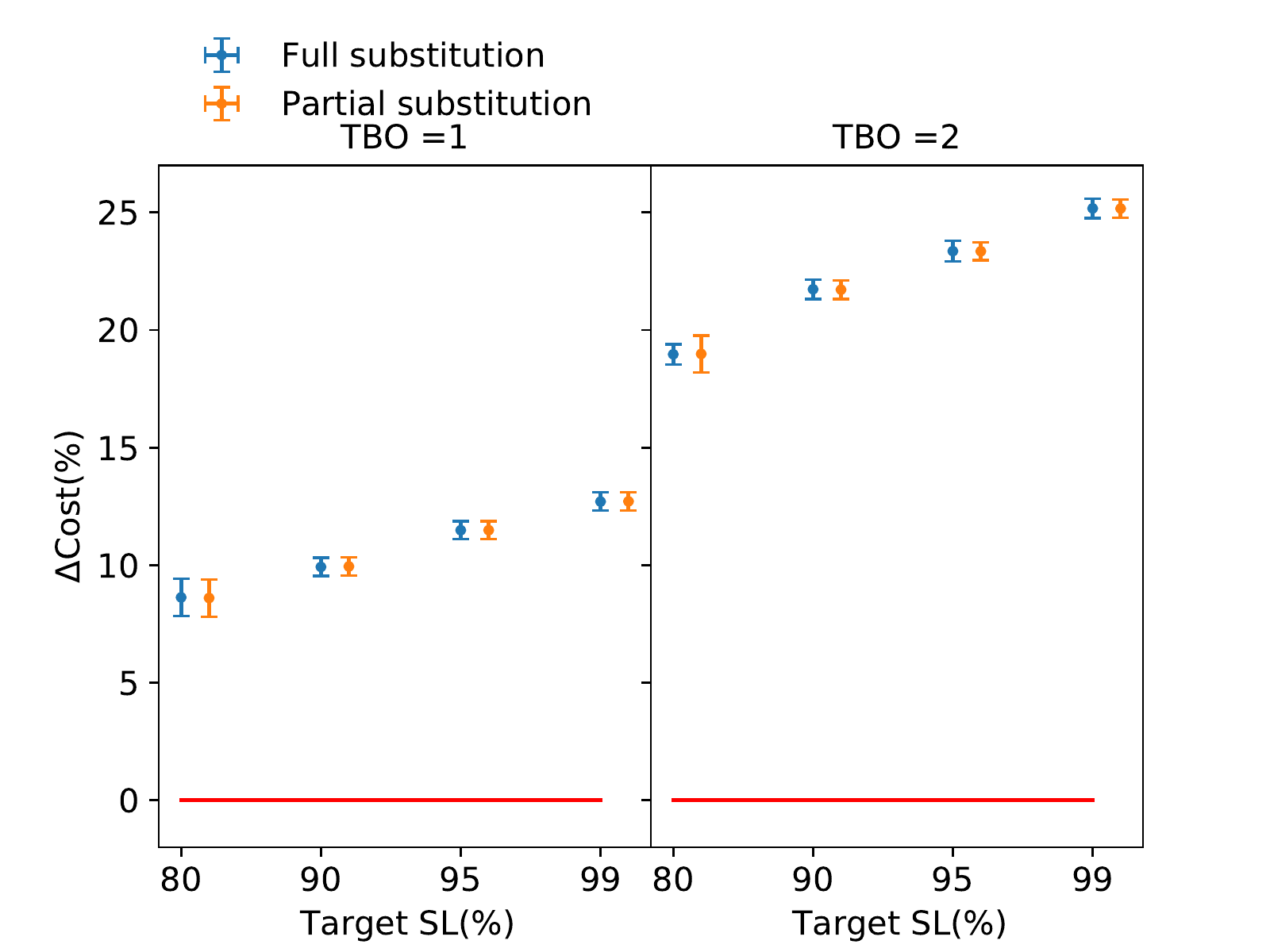}
\caption{Effect of substitution (Relative cost decrease)} 
\label{fig:sustitution}
\end{center}
\end{figure}

In Figure~\ref{fig:sustitution} we display the cost reduction of the two cases in which substitution is allowed relative to the case with no substitution, for varying values of service level and TBO equal to 1 and 2. These results indicate that substitution enables significant cost savings, and that the savings are significantly higher when TBO is higher (i.e., for instances where the setup costs are higher).  We also observe that the cost savings are about the same with full and partial substitution, showing that a limited amount of allowed product substitution can capture the majority of the benefit. 




\section{Conclusion}
\label{sec:conc}

We study an infinite-horizon stochastic lot-sizing problem with a supplier-driven product substitution option and the service level constraint which is defined jointly over different products. To solve this problem, we consider a finite-horizon version of this problem and apply it in a rolling-horizon framework. We propose different policies based on solving a different approximation of this multi-stage problem to make the decisions in each period.  We propose two deterministic policies and a policy based on solving a two-stage chance-constrained stochastic program. We also present a branch-and-cut algorithm for effectively solving the two-stage chance-constrained model. 

We conducted an extensive evaluation comparing these policies within a simulation study. The results indicate that the proposed chance-constraint policy leads to reliable satisfaction of the service level targets, and does so at significantly lower cost than the approximations based on solving deterministic models. Most significantly, we find that allowing supplier-driven substitution can lead to very significant reductions in costs to meet a desired service level target, and that these reductions can be obtained by allowing product substitution between a relatively limited range of products, and are most significant when setup costs are relatively higher.

\ACKNOWLEDGMENT{
The authors gratefully acknowledge the support of the Digital Research Alliance of Canada and FRQNT International Internship Program.
}

\bibliographystyle{informs2014} 
\bibliography{References}

\newpage
\bigskip
\begin{APPENDICES}
\small
\section{Cut Separation Algorithm}
\label{APP1}
\begin{algorithm}[ht]
\SetAlgoLined
\linespread{1.15}\selectfont
{OUTPUT: A most violated mixing inequality defined by the ordered index set $T$ } \\
INPUT: $\zsol_{\omega}$, 
$h_{\omega}(\pisol,\betasol)$ for $\omega \in \Omega$, $p$
\\
Sort the $\zsol$ components to obtain permutation $\sigma$ of the indices satisfying:
$\zsol_{\sigma_1} \leq \zsol_{\sigma_2} \leq \cdots \leq \zsol_{\sigma_{p+1}} \leq \cdots$ \\
$v \gets h_{\sigma_{p+1}}(\pisol,\betasol)$ \\
$T \gets \{ \}$ \\
$ i \gets 1$ \\
\While{$v < h_{\sigma_1}(\pisol,\betasol)$}{
  \If{$h_{\sigma_i}(\pisol,\betasol) > v$}{
  $T \gets T \cup \{\sigma_i\}$ \\
  $v \gets h_{\sigma_i}(\pisol,\betasol)$  }
  $i \gets i+1$
}
\caption{Finding a most violated inequality of the form \eqref{eq:mixed}.}\label{alg:mostviolated}
\end{algorithm}
\end{APPENDICES}

\end{document}

%% file: Definitions.tex
\usepackage{epsfig}
\usepackage{lscape}
\usepackage{adjustbox}
\usepackage{mathtools}
\usepackage{amsmath}
\usepackage{comment}
\usepackage{booktabs}
\usepackage{graphicx}
\usepackage{array}
\usepackage{xcolor}
\usepackage[margin=1 in]{geometry}
\usepackage[ruled,vlined,noend]{algorithm2e}
\usepackage{amssymb}
\usepackage[colorlinks=true,linkcolor=blue,citecolor=blue, hypertexnames=false]{hyperref} 
\usepackage{textcomp}
\usepackage{caption}
\usepackage{multirow}
\usepackage{here}
\usepackage{amsmath}
\usepackage{dsfont}

\usepackage{subfig}
\allowdisplaybreaks
 
\newcommand{\ti}{t} 
\newcommand{\TI}{\mathcal{T}}
\newcommand{\Ti}{T}
\newcommand{\ka}{k} 
\newcommand{\KA}{\mathcal{K}}
\newcommand{\Ka}{K}
\newcommand{\jey}{j} 
\newcommand{\Bi}{b} 
\newcommand{\Vi}{v} 
\newcommand{\Es}{s}
\newcommand{\Zed}{z} 
\newcommand{\m}{\omega} 
\newcommand{\Em}{|\Omega|} 
\newcommand{\EM}{\Omega} 

\newcommand{\x}{x} 
\newcommand{\y}{y} 
\newcommand{\inv}{i} 
\newcommand{\InvPos}{available inventory}

\newcommand{\Csub}{\mathcal{K}^+_k}
\newcommand{\Psub}{\mathcal{K}^-_k}

\newcommand{\tHat}{\hat{\ti}} 
\newcommand{\tAct}{\hat{\ti}}
\newcommand{\PSpolicy}{``production substitution policy"}

\newcommand{\Sc}{\overline{s}}
\newcommand{\Bc}{\overline{b}}

\newcommand{\cred}{\color{black}}

%% file: LiteratureTable.tex
\begin{table}[ht]
\caption{Summary of related papers.}
\begin{adjustbox}{width=1\textwidth,center=\textwidth}
\begin{tabular}{llllllllll}
\toprule
 & & Planning \\
                          & Year                                    & Horizon                        & Uncertainty                       & Problem                         & Subst.        & Service level                & Strategy       & Model               &      Method.        \\
                          \midrule

\citeauthor{bitran1992deterministic}         & 1992                                    & F                                       & Yield                             & Co-production                   & SD                  & J (Products)                 & S,D            & LP                          &     A         \\
\citeauthor{bitran1992ordering}           & 1992                                    & I                                       & Yield                             & Co-production                   & SD                  & I                            & D              & LP                      &  H            \\
\citeauthor{bassok1999single}             & 1999                                    & S                                       & Dem, Yield                     & Periodic review inventory & SD                  &                              &                &                           &       G       \\
\citeauthor{hsu1999random}            & 1999                                    & S                                       & Dem, Yield                     & Co-production                   & SD                  &                              &                & MILP                   &    G           \\
\citeauthor{rao2004multi}                & 2004                                    &   S                                      & Dem                               & Inventory planning + setup      & SD                  &                              &                &                 MILP         &  H            \\
\citeauthor{hsu2005dynamic}                & 2005                                    & F                                       & Det                               & Lot sizing                      & SD                  &                              &                &                  MILP         & DP             \\
\citeauthor{nagarajan2008inventory} & 2008                                    & S, F                                 & Dem                               & Inventory planning              &CD                  &                              & D              &                         &  H             \\
\citeauthor{lang2010efficient}         & 2010                                    & F                                       & Det                               & Lot sizing                      & SD                  &                              &                & MILP                        &              \\
\citeauthor{ng2012robust}              & 2012                                    & S                                       & Dem                               & Co-production                   &CD                  & \multicolumn{2}{l}{M}  &                    LP, MILP       &              \\
\citeauthor{zhang2014branch}             & 2014                                    & F                                       & Dem                               & Inventory planning                      &                     & J (Periods)                   & D              &              MILP        & B\&C               \\
\citeauthor{jiang2017production}              & 2017                                    & F                                       & Dem                               & \multicolumn{2}{l}{Production planning}               & J (Periods)                   & S              &      MILP                  &  SAA            \\
\citeauthor{gicquel2018joint}         & 2018                                    & F, I                                  & Dem                               & Lot sizing                      &                     & J (Periods)                   & S              &  MILP                 &   SAA           \\
\citeauthor{liu2018polyhedral}          & 2018                                    & F                                       & Dem                               & Lot sizing                      &                     & J (Periods)                   & S              &                 MILP   & B\&C                \\
\citeauthor{chen2020dynamic}             & 2020                                    &    F                                     & Dem                               & Inventory control               &CD                  &                              &                & &  OL \\
\citeauthor{akccaycategory}               & 2020                                    & S                                       & Dem                               & Inventory planning                                 &CD                  & I                            &                & & A\\
Our work                  &                                         & I                                       & Dem                               & Lot sizing                      & SD                  & J (Products)                 & D              & MILP                &           B\&C   

\\
\midrule
\textbf{Acronyms}                  &                                         &                                         &                                   &                                 &                     &                              &                &                           &              \\
\multicolumn{10}{l}{Planning Horizon ..          I: infinite,                              F: Finite,                               S: Single period}                                           \\
\multicolumn{10}{l}{
Uncertainty ..               Det: Deterministic,                       Dem: Random Demand,                      Yield : Random Yield}\\

\multicolumn{10}{l}{
Substitution ..              SD: Supplier-driven,                          CD: Customer-driven }                                     \\
\multicolumn{10}{l}{
Service level  ..            I: Individual,                            J (Periods): Joint over multiple periods,  J (Products): Joint over multiple   products,  M: Maximizing service level  }         \\
\multicolumn{10}{l}{
Strategy    ..               S: Static,                                D: Dynamic  }               \\
\multicolumn{10}{l}{
Model   ..             MILP: Mixed-integer linear programming,  LP: Linear programming }             \\

\multicolumn{10}{l}{
 Methodology .. B\&C: Branch-and-cut algorithm,             G: Greedy algorithm,  A: Model approximation,    H: Heuristics, DP: Dynamic programming,                         } \\
\multicolumn{10}{l}{ SAA: Sample average approximation, OL: Online learning  }               \\
  \bottomrule
\end{tabular}
\end{adjustbox}
\label{tab:litRev}
\end{table}